\documentclass[11pt,reqno]{amsart}
\usepackage{amssymb,amscd,amsbsy}
\usepackage{amssymb,amscd,amsbsy,mathrsfs}
\newcommand{\lb}{\linebreak}

\renewcommand{\a}{\alpha}

\renewcommand{\d}{\delta}
\newcommand{\e}{\varepsilon}

\newcommand{\z}{\zeta}

\renewcommand{\l}{\lambda}

\newcommand{\s}{\sigma}
\renewcommand{\t}{\tau}

\newcommand{\f}{\varphi}
\renewcommand{\o}{\omega}

\newcommand{\G}{\Gamma}
\newcommand{\D}{\Delta}
\renewcommand{\L}{\Lambda}

\renewcommand{\O}{\Omega}

\newcommand{\A}{{\mathscr A}}

\newcommand{\cd}{{\mathscr D}}
\newcommand{\F}{{\mathscr F}}

\newcommand{\h}{{\mathscr H}}

\newcommand{\N}{{\mathscr N}}

\newcommand{\cR}{{\mathscr R}}

\newcommand{\X}{{\mathscr X}}
\newcommand{\Y}{{\mathscr Y}}
\newcommand{\cU}{{\mathscr U}}

\newcommand{\C}{{\Bbb C}}

\newcommand{\dd}{{\Bbb D}}
\newcommand{\R}{{\Bbb R}}
\newcommand{\Z}{{\Bbb Z}}

\newcommand{\0}{{\boldsymbol{0}}}

\newcommand{\bs}{\boldsymbol}

\newcommand{\m}{{\boldsymbol m}}
\newcommand{\bS}{{\boldsymbol S}}

\newcommand{\rf}[1]{(\ref{#1})}

\newcommand{\df}{\stackrel{\mathrm{def}}{=}}
\newcommand{\dist}{\operatorname{dist}}

\newcommand{\re}{\operatorname{Re}}

\newcommand{\supp}{\operatorname{supp}}
\newcommand{\clos}{\operatorname{clos}}

\newcommand{\rank}{\operatorname{rank}}
\newcommand{\const}{\operatorname{const}}

\newcommand{\eeq}{\end{equation}}
\newcommand{\beq}{\begin{equation}}
\newcommand{\bay}{\begin{eqnarray}}
\newcommand{\ba}{\begin{align*}}
\newcommand{\ea}{\end{align*}}
\newcommand{\ey}{\end{eqnarray}}
\newcommand{\bey}{\begin{eqnarray*}}
\newcommand{\eey}{\end{eqnarray*}}

\newcommand{\imp}{\Rightarrow}
\newcommand{\be}{\infty}

\newcommand{\bl}{\blacksquare}

\newcommand{\Pf}{{\bf Proof. }}
\newcommand{\im}{\operatorname{Im}}
\renewcommand{\re}{\operatorname{Re}}
\newcommand{\ov}{\overline}

\newtheorem{thm}{\hspace{\parindent}Theorem}[section]

\newtheorem{cor}[thm]{\hspace{\parindent}Corollary}
\newtheorem{lem}[thm]{\hspace{\parindent}Lemma}

\theoremstyle{remark}

\newtheorem*{rem*}{Remark}

\newcommand\Li{{\rm Lip}}
\newcommand\fM{\frak M}

\newcommand\dg{\frak D}

\newcommand\mX{\mathcal{X}}
\newcommand\mY{\mathcal{Y}}
\newcommand\mB{\mathcal{B}}
\newcommand\fn{\frak n}
\newcommand\fm{\frak m}
\newcommand\fp{\frak p}
\newcommand\mZ{\mathcal Z}

\newcommand\mI{\mathcal{I}}

\newcommand{\OL}{{\rm OL}}
\newcommand{\COL}{{\rm CL}}

\newcommand{\SA}{{\bf SA}}
\newcommand{\No}{{\bf N}}
\newcommand{\Un}{{\bf U}}
\newcommand{\Pro}{{\bf P}}
\newcommand{\USA}{{\bf USA}}
\newcommand{\fF}{{\frak F}}
\newcommand{\CO}{{\bf C}}
\newcommand{\sgn}{\operatorname{sgn}}
\newcommand{\card}{\operatorname{card}}

\begin{document}

\newcommand{\vse}{\vspace{.2in}}
\numberwithin{equation}{section}

\title{Operator and commutator moduli\\ of continuity for normal operators}
\author{A.B. Aleksandrov and V.V. Peller}

\begin{abstract}
We study in this paper properties of functions of perturbed normal operators
and develop earlier results obtained in \cite{APPS2}. We study operator Lipschitz and
commutator Lipschitz functions on closed subsets of the plane. For such functions
we introduce the notions of the operator modulus of continuity and of various
commutator moduli of continuity. Our estimates lead to estimates of the norms
of quasicommutators $f(N_1)R-Rf(N_2)$ in terms of $\|N_1R- RN_2\|$, where $N_1$ and
$N_2$ are normal operator and $R$ is a bounded linear operator. In particular, we show
that if $0<\a<1$ and $f$ is a H\"older function of order $\a$, then for normal operators
$N_1$ and $N_2$,
$$
\|f(N_1)R-Rf(N_2)\|\le\const(1-\a)^{-2}\|f\|_{\L_\a}\|N_1R-RN_2\|^\a\|R\|^{1-\a}.
$$
In the last section we obtain lower estimates for constants in operator H\"older
estimates.
\end{abstract}

\maketitle

\

\begin{center}
{\Large Contents}
\end{center}

\

\begin{enumerate}
\item[1.] Introduction \quad\dotfill \pageref{In}
\item[2.] Schur multipliers and double operator integrals \quad\dotfill \pageref{Schm}
\item[3.] Operator Lipschitz and commutator Lipschitz functions
\quad\dotfill \pageref{OLCL}
\item[4.] Absolutely convergent Fourier integrals
and estimates of commutator Lipschitz norms  \quad\dotfill \pageref{FouIn}
\item[5.] Operator and commutator moduli of continuity \quad\dotfill \pageref{Omc}
\item[6.] Estimates of commutator moduli of continuity \quad\dotfill \pageref{ECMC}
\item[7.] Constants in operator H\"older inequalities \quad\dotfill \pageref{alpha}
\item[] References \quad\dotfill \pageref{bibl}
\end{enumerate}

\

\setcounter{section}{0}
\section{\bf Introduction}
\setcounter{equation}{0}
\label{In}

\medskip

In this paper we continue the study of properties of functions of normal operators under perturbation. This study was undertaken in
\cite{APPS2} (see also \cite{APPS}). Let us summarize briefly some results obtained in \cite{APPS2}.

It was shown in \cite{APPS2} that if $f$ is a function on the plane that belongs to the (homogeneous) Besov space $B_{\be1}^1(\R^2)$, then it is {\it operator Lipschitz}, i.e.,
$$
\|f(N_1)-f(N_2)\|\le\const\|N_1-N_2\|
$$
for arbitrary normal (not necessarily bounded) operators $N_1$ and $N_2$ such that $N_1-N_2$ is bounded. We refer the reader to \cite{Pee} for definitions and basic properties of Besov spaces.

Note that a {\it Lipschitz function} on the plane (i.e., a function $f$ such that
$$
|f(\z_1)-f(\z_2)|\le\const|\z_1-\z_2|,\quad\z_1,\,\z_2\in\C)
$$
does not have to be operator Lipschitz. This is not true even for functions defined on the real line $\R$, i.e., there are Lipschitz functions $f$ on $\R$ such that
$$
\sup\frac{\|f(A)-f(B)\|}{\|A-B\|}=\be,
$$
where the supremum is taken over bounded self-adjoint operators $A$ and $B$.
The first example of such functions was constructed in \cite{F}.
Later it was shown in \cite{Mc} and \cite{K} that the function $x\mapsto|x|$ on $\R$ is not
operator Lipschitz. Note also that in \cite{Pe2} a necessary condition was found: if
a function $f$ on $\R$ is operator Lipschitz, then $f$ belongs locally to the Besov class
$B_{1,1}^1(\R)$.

It was also shown in \cite{APPS2} that if $f$ belongs to the {\it H\"older class} $\L_\a(\R^2)$, $0<\a<1$, i.e.,
$$
|f(\z_1)-f(\z_2)|\le\const|\z_1-\z_2|^\a,\quad\z_1,\,\z_2\in\C,
$$
then $f$ is {\it operator H\"older of order} $\a$, i.e.,
\bay
\label{a^-1}
\|f(N_1)-f(N_2)\|\le c_\a\|f\|_{\L_\a}\|N_1-N_2\|^\a
\ey
for arbitrary (not necessarily bounded) normal operators $N_1$ and
$N_2$ such that $N_1-N_2$ is bounded. Moreover, it is shown in \cite{APPS2} that
$c_\a\le\const(1-\a)^{-1}$. Here
$$
\|f\|_{\L_\a}\df\sup_{\z_1\ne\z_2}\frac{|f(\z_1)-f(\z_2)|}{|\z_1-\z_2|^\a}.
$$

More general results were also obtained in \cite{APPS2} in
the case of arbitrary moduli of continuity.
Recall that a continuous nondecreasing function
$\o:[0,\be)\to[0,\be)$ is called a {\it modulus of continuity} if
$$
\o(x+y)\le\o(x)+\o(y),\quad x,\,y\in[0,\be).
$$
It was shown in \cite{APPS2} that if $\o$ is a modulus of continuity and $f\in\L_\o(\R^2)$, i.e.,
$$
|f(\z_1)-f(\z_2)|\le\const\o\big(|\z_1-\z_2|\big),\quad\z_1,\,\z_2\in\C,
$$
then
$$
\|f(N_1)-f(N_2)\|\le C\|f\|_{\L_\o}\,\o_*\big(\|N_1-N_2\|\big),
$$
where $C$ is a numerical constant and
\bay
\label{o*}
\o_*(\d)\df\d\int_\d^\be\frac{\o(t)}{t^2}\,dt=
\int_1^\be\frac{\o(t\d)}{t^2}\,dt,
\quad\d>0.
\ey
Here
$$
\|f\|_{\L_\o}\df\sup_{\z_1\ne\z_2}\frac{|f(\z_1)-f(\z_2)|}{\o(|\z_1-\z_2|)}.
$$

We would like to also mention here that in \cite{APPS2} estimates for Schatten--von Neumann norms
norms of $f(N_1)-f(N_2)$ as well as other ideal norms are also obtained.

Note that similar results were obtained earlier for self-adjoint operators
(as well as for unitary operators, contractions, and dissipative operators)
in \cite{Pe2}, \cite{Pe3}, \cite{Pe5}, \cite{AP1}, \cite{AP2}, \cite{AP3},
\cite{AP4}, \cite{Pe0}, \cite{AP5}.

Analogous estimates were obtained in \cite{APPS2} for commutators and
quasicommutators. Namely, it was shown in \cite{APPS2} that if $f\in
B_{\be1}^1(\R^2)$, then \bay \label{bes}
\|f(N_1)R-Rf(N_2)\|\le\const\|f\|_{B_{\be1}^1}
\max\big\{\|N_1R-RN_2\|,\|N^*_1R-RN^*_2\|\big\} \ey for an arbitrary
bounded operator $R$ and arbitrary normal operators $N_1$ and $N_2$
such that the operators $N_1R-RN_2$ and $N^*_1R-RN^*_2$ are bounded.

Next, the following analog of \rf{a^-1} was proved in \cite{APPS2}:
\bay \label{chol} \|f(N_1)R-Rf(N_2)\| \le c_\a\|f\|_{\L_\a}
\max\big\{\|N_1R-RN_2\|,\|N^*_1R-RN^*_2\|\big\}^\a\|R\|^{1-\a}, \ey
whenever $f\in\L_\a(\R^2)$, $R$ is a bounded operator, and $N_1$ and
$N_2$ are normal operators such that the operators $N_1R-RN_2$ and
$N^*_1R-RN^*_2$ are bounded. Here $c_\a\le\const(1-\a)^{-1}$.

Finally, it was shown in \cite{APPS2} that under the same hypotheses on $R$, $N_1$ and $N_2$,
the following inequality holds for an arbitrary modulus of continuity $\o$ and
an arbitrary function $f$ in $\L_\o(\R^2)$:
\begin{align}
\label{com}
\|f(N_1)R&-Rf(N_2)\|\\[.2cm]
& \le\const\|f\|_{\L_\o} \|R\|\,\o_*\!
\left(\frac{\max\big\{\|N_1R-RN_2\|,\|N^*_1R-RN^*_2\|\big\}}{\|R\|}\right).\nonumber
\end{align}

In this paper we consider the problem of whether we can estimate the
quasicommutator norms $\|f(N_1)R-Rf(N_2)\|$ in terms of
$\|N_1R-RN_2\|$ rather than in terms of
\lb$\max\big\{\|N_1R-RN_2\|,\|N^*_1R-RN^*_2\|\big\}$.

Let us first mention that in inequality \rf{bes} it is impossible to
replace
\lb$\max\big\{\|N_1R-RN_2\|,\|N^*_1R-RN^*_2\|\big\}$ with
$\|N_1R-RN_2\|$. Indeed, it can be deduced from results of \cite{JW}
that if
\bay
\label{krysa}
\|f(N_1)R-Rf(N_2)\|\le \const \|N_1R-RN_2\|
\ey
for arbitrary bounded $N_1$ and $N_2$ with spectra contained in a
given closed set ${\frak F}$, then $f$ must have complex
derivative at each nonisolated point of $\fF$. In particular, if $\fF=\C$ and
$f$ satisfies \rf{krysa}, then $f(z)=az+b$ for some $a,\,b\in\C$.

Surprisingly, it turns out that inequality \rf{chol} still holds if
we replace \lb$\max\big\{\|N_1R-RN_2\|,\|N^*_1R-RN^*_2\|\big\}$ with
$\|N_1R-RN_2\|$. However, the constant $c_\a$ jumps. Namely, we show
in \S\,\ref{ECMC} of this paper that
\bay
\label{bezm}
 \|f(N_1)R-Rf(N_2)\|
\le c_\a\|f\|_{\L_\a} \|N_1R-RN_2\|^\a\|R\|^{1-\a}
\ey
with
$c_\a\le\const(1-\a)^{-2}$. We do not know whether inequality \rf{bezm} holds with
$c_\a\le\const(1-\a)^{-1}$.

We also obtain in \S\,\ref{ECMC} the following modification of inequality \rf{com}:
\bay
\label{woutm}
\|f(N_1)R-Rf(N_2)\|
\le\const\|f\|_{\L_\o} \|R\|\,\o_{**}\!
\left(\frac{\|N_1R-RN_2\|}{\|R\|}\right),
\ey
where $\o_{**}\df(\o_*)_*$. Again, we do not know whether we can replace in \rf{woutm}
$\o_{**}$ with $\o_*$.

In \S\,\ref{alpha} we study the problem of whether
our estimate of the constant $c_\a$ in
inequality \rf{a^-1} is sharp. We show that $c_\a\ge C(1-\a)^{-1/2}$
for a positive number $C$.

We introduce in \S\,\ref{Omc} various commutator moduli
of continuity and study their properties. We
study in \S\,\ref{OLCL} some properties of operator
Lipschitz and commutator Lipschitz
functions.

In \S\,\ref{FouIn} we give some auxiliary results: norm estimates in the space of
functions with absolutely convergent Fourier integrals and estimates of commutator
Lipschitz norms.

Finally, in \S\,\ref{Schm} we give an introduction into Schur multipliers and double
operator integrals.

\

\section{\bf Schur multipliers and double operator integrals}
\setcounter{equation}{0}
\label{Schm}

\

We define in this section notion of Schur multipliers associated with
two spectral measures. However, we start the section with the definition
of Schur multipliers associated with two scalar measures.
This corresponds to the case of spectral measures of multiplicity 1.
We discuss properties of Schur multipliers and define double operator
integrals.

Let $(\mathcal X,\mu)$ and $(\mathcal Y,\nu)$ be $\s$-finite measure spaces.
Let $k\in L^2(\mathcal X\times\mathcal Y,\mu\otimes\nu)$. Then $k$ induces
the integral operator $\mathcal I_k=\mathcal I_k^{\mu,\nu}$ from
$L^2(\nu)=L^2(\mathcal Y,\nu)$ to $L^2(\mu)=L^2(\mathcal X,\mu)$ defined by
$$
(\mathcal I_k f)(x)=\int_{\mathcal Y}k(x,y)f(y)\,d\nu(y),\quad f\in L^2(\mathcal Y,\nu).
$$
Denote by $\|k\|_{\mathcal B}=\|k\|_{\mB_{\mX,\mY}^{\,\mu,\nu}}$ the operator norm
of $\mathcal I_k$.
Let $\Phi$ be a $\mu\otimes\nu$-measurable function defined almost everywhere on
$\mX\times\mY$.
We say that $\Phi$ is a {\it Schur multiplier with respect to $\mu$ and $\nu$} if
$$
\|\Phi\|_{\frak M_{\mX,\mY}^{\,\mu,\nu}}
\df\sup\big\{\|\Phi k\|_{\mB}:
~k\,,\Phi k\in L^2(\mathcal X\times\mathcal Y,\mu\otimes\nu),~\|k\|_{\mB}\le1\big\}<\be.
$$
We denote by $\frak M_{\mX,\mY}^{\,\mu,\nu}$ the space of Schur multipliers with respect to $\mu$ and $\nu$.
It can be shown easily that $\frak M_{\mX,\mY}^{\,\mu,\nu}\subset L^\be(\mX\times\mY,\mu\otimes\nu)$ and
\bay
\label{bemult}
\|\Phi\|_{L^\be(\mX\times\mY,\mu\otimes\nu)}\le
\|\Phi\|_{\frak M_{\mX,\mY}^{\,\mu,\nu}}.
\ey
%In particular, in the definition of $\|\cdot\|_{\frak M_{\mX,\mY}^{\,\mu,\nu}}$ we can omit the condition $\Phi k\in L^2(\mathcal X\times\mathcal Y,\mu\otimes\nu)$.
Thus if $\Phi\in\frak M_{\mX,\mY}^{\,\mu,\nu}$, then
$$
\|\Phi\|_{\frak M_{\mX,\mY}^{\,\mu,\nu}}
=\sup\big\{\|\Phi k\|_{\mB}:
~k\in L^2(\mathcal X\times\mathcal Y,\mu\otimes\nu),~\|k\|_{\mB}\le1\big\}.
$$
Note that $\frak M_{\mX,\mY}^{\,\mu,\nu}$ is a Banach algebra:
$$
\|\Phi_1\Phi_2\|_{\frak M_{\mX,\mY}^{\,\mu,\nu}}
\le\|\Phi_1\|_{\frak M_{\mX,\mY}^{\,\mu,\nu}}\|\Phi_2\|_{\frak M_{\mX,\mY}^{\,\mu,\nu}}.
$$
It is easy to see that $\|\Phi\|_{\frak M_{\mX,\mY}^{\,\mu,\nu}}=\|\Phi_\flat\|_{\frak M_{\mY,\mX}^{\,\nu,\mu}}$
for $\Phi_\flat(y,x)=\Phi(x,y)$.

Note that if $\mX$ and $\mY$ are at most countable sets, and $\mu$ and $\nu$
are the counting measures on $\mX$ and $\mY$, the above definition coincides
with the definition of Schur multipliers on the space of matrices:
a matrix $A=\{a_{jk}\}$ is called a Schur multiplier on the space of
bounded matrices if
$$
A\star B\quad\mbox{is a matrix of a bounded operator, whenever}\quad B\quad\mbox{is}.
$$
Here we use the notation
\bay
\label{ScHad}
A\star B=\{a_{jk}b_{jk}\}
\ey
for the Schur--Hadamard product of the matrices $A=\{a_{jk}\}$ and $B=\{b_{jk}\}$

Clearly, the norm of $A$ in the space of Schur multipliers is the norm of
the transformer
$$
B\mapsto A\star B.
$$

We need the following known result:

\begin{lem}
\label{fyokla}
Let $\{G_n\}_{n=1}^\be$ be a sequence of disjoint
$\mu$-measurable subsets of $\mX$ and let $\{H_n\}_{n=1}^\be$ be
a sequence of disjoint $\nu$-measurable subsets
of $\mY$. Put ${Z}\df\bigcup\limits_{n=1}^\be (G_n\times H_n)$.
Then $\|\chi_{_{{Z}}}\|_{\frak M_{\mX,\mY}^{\,\mu,\nu}}\le1$.
\end{lem}

\Pf
Let $k\in L^2(\mX\times\mY,\mu\otimes\nu)$, $f\in L^2(\mY,\nu)$ and
$g\in L^2(\mX,\mu)$. We have
\begin{align*}
|(\mI_{\chi_{_{{Z}}}k}f,g)|&=
\left|\sum_{n=1}^\be\big(\mI_k (\chi_{_{H_n}}f),\chi_{_{G_n}}g\big)\right|\\[.2cm]
&\le\|k\|_{\mB}\sum_{n=1}^\be\|\chi_{_{H_n}}f\|_{L^2(\nu)}\|\chi_{_{G_n}}f\|_{L^2(\mu)}\\
&\le\|k\|_{\mB}\left(\sum_{n=1}^\be\|\chi_{_{H_n}}f\|_{L^2(\nu)}^2\right)^{1/2}
\left(\sum_{n=1}^\be\|\chi_{_{G_n}}g\|_{L^2(\mu)}^2\right)^{1/2}\\[.2cm]
&\le
\|k\|_{\mB}\|f\|_{L^2(\nu)}\|g\|_{L^2(\mu)}.
\end{align*}
Hence, $\|\chi_{_{{Z}}}k\|_{\mB}\le\|k\|_{\mB}$, and so
$\|\chi_{_{{Z}}}\|_{\frak M_{\mX,\mY}^{\,\mu,\nu}}\le1$. $\bl$

Clearly, taking \rf{bemult} into account, we find that
\bay
\label{1mult}
\|\chi_{_{{Z}}}\|_{\frak M_{\mX,\mY}^{\,\mu,\nu}}=1,
\ey
whenever
$(\mu\otimes\nu)({Z})=\sum_{n=1}^\be\mu(G_n)\nu(H_n)>0$.

To state a description of the space of Schur multipliers we define
the {\it integral projective tensor product} $L^\be(\mu)\hat\otimes_{\rm i}
L^\be(\nu)$ of $L^\be(\mu)$ and $L^\be(\nu)$. We say that
$\Phi\in L^\be(\mu)\hat\otimes_{\rm i}
L^\be(\nu)$ if $\Phi$ admits a representation
\bay
\label{ipt}
\Phi(x,y)=\int_\O \f(x,w)\psi(y,w)\,d\l(w),
\ey
where $(\O,\l)$ is a $\s$-finite measure space, $\f$ is a measurable function on $\X\times \O$,
$\psi$ is a measurable function on $\Y\times \O$, and
$$
\int_\O\|\f(\cdot,w)\|_{L^\be(\mu)}\|\psi(\cdot,w)\|_{L^\be(\nu)}\,d\l(w)<\be.
$$

The space of Schur multipliers admits the following description:

\medskip

{\bf Theorem on Schur multipliers.} {\em Let $(\mX,\mu)$ and $(\mY,\nu)$
be $\s$-finite measure spaces and let
$\Phi$ be a measurable function on
$\X\times\Y$. The following are equivalent:

{\rm (i)} $\Phi\in\frak M_{\mX,\mY}^{\,\mu,\nu}$;

{\rm (ii)} $\Phi\in L^\be(\mu)\hat\otimes_{\rm i}L^\be(\nu)$;

{\rm (iii)} there exist a $\s$-finite measure space $(\O,\l)$,
 measurable functions $\f$ on $\X\times\O$ and $\psi$ on $\Y\times\O$ such that
{\em\rf{ipt}} holds and
$$
\left\|\left(\int_\O|\f(\cdot,w)|^2\,d\l(w)\right)^{1/2}\right\|_{L^\be(E)}
\left\|\left(\int_\O|\psi(\cdot,w)|^2\,d\l(w)\right)^{1/2}\right\|_{L^\be(F)}<\be.
$$
}

We refer the reader to \cite{Pe2} for more detailed information and references.

Let $\mX$ and $\mY$ be closed subsets of $\C$.
We denote by ${\frak M}_{\mX,\mY}$ the space of
Borel Schur multipliers on $\mX\times\mY$, i.e., the space of
Borel functions $\Phi$ defined everywhere on $\mX\times\mY$ such that
$$
\|\Phi\|_{\frak M_{\mX,\mY}}\df\sup\|\Phi\|_{\frak M_{\mX,\mY}^{\mu,\nu}}<\be,
$$
where the supremum is taken over all Borel measures
$\mu$ and $\nu$ on $\mathcal X$ and $\mathcal Y$. In the case  $\mX=\mY$, we use
the notation
$$
{\frak M}_\mX \df{\frak M}_{\mX,\mX}.
$$

It can be shown easily that
$$
\sup_{(x,y)\in\mX\times\mY}|\Phi(x,y)|\le\|\Phi\|_{\frak M_{\mX,\mY}}.
$$
It is also easy to verify that if $\Phi_n\in\fM_{\mX,\mY}$,
$\Phi$ is a bounded Borel function on $\mX\times\mY$, and
$\Phi_n(x,y)\to \Phi(x,y)$ for all $(x,y)\in\mathcal X\times\mathcal Y$,
then
$$
\|\Phi\|_{\frak M_{\mX,\mY}}
\le\liminf\limits_{n\to\be}
\|\Phi_n\|_{\frak M_{\mX,\mY}}.
$$
In particular, $\Phi\in\frak M_{\mX,\mY}$ if
$\liminf\limits_{n\to\be}
\|\Phi_n\|_{\frak M_{\mX,\mY}}<\be$.

We need the following version of Lemma \ref{fyokla}.

\begin{lem}
Let $\mu$ and $\nu$ be Borel measures on closed subsets $\mX$ and $\mY$ of $\C$.
Put
$$
\D\df\{(x,y)\in\mX\times\mY~:x=y\}.
$$
Then
$\|\chi_{_{{\D}}}\|_{\frak M_{\mX,\mY}^{\,\mu,\nu}}\le1$
and hence, $\|\chi_{_{{\D}}}\|_{\frak M_{\mX,\mY}}\le1$.
\end{lem}

\Pf Put $\mX_0\df\{x\in\mX:\mu(\{x\})>0\}$,
$\mY_0\df\{y\in\mY:\nu(\{y\})>0\}$
and $\D_0\df\{(x,y)\in\mX_0\times\mY_0:x=y\}$. Clearly, $\mX_0$ and $\mY_0$
are at most countable. It is easy to see that
$(\mu\otimes\nu)(\D\setminus\D_0)=0$.
Hence, $\|\chi_{_{{\D}}}\|_{\frak M_{\mX,\mY}^{\,\mu,\nu}}=
\|\chi_{_{{\D_0}}}\|_{\frak M_{\mX,\mY}^{\,\mu,\nu}}$. It remains to observe
that $\|\chi_{_{{\D_0}}}\|_{\frak M_{\mX,\mY}^{\,\mu,\nu}}\le1$
by Lemma \ref{fyokla}. $\bl$

\begin{cor}
\label{grunya}
Under the hypotheses of the lemma, the following inequality holds:
$$
\|\chi_{_{{(\mX\times\mY)\setminus\D}}}\|_{\fM_{\mX,\mY}^{\,\mu,\nu}}\le2
\quad\text{and}\quad\|\chi_{_{{(\mX\times\mY)\setminus\D}}}\|_{\fM_{\mX,\mY}}\le2.
$$
\end{cor}

%<<<<<<<<<<<<<<<<<<<<<
%We are going to deal with functions $f$ on $\mX\times\mY$ that are continuous in each variable.
%It must be a well-known fact that such a function $f$ has to be a Borel function.
%Indeed, one can construct an increasing sequence $\{\mY_n\}_{n=1}^\be$
%of discrete closed subsets of $\mY$ such that $\bigcup\limits_{n=1}^\be\mY_n$
%is dense in $\mY$. Let us consider the function $f_n:\mX\times\R\to\C$
%such that $f\big|(\mX\times\mY_n)=f_n\big|(\mX\times\mY_n)$ and $f_n(x,\cdot)$ is
%a piecewise linear function with nodes in $\mY_n$ for all $x\in\mX$.
%Clearly, the function $f_n$ is defined uniquely if we require
%that $f_n(x,\cdot)$ is constant on each unbounded
%complimentary interval of $\mY_n$.
%It is easy to see that $f_n$ is continuous on $\mX\times\R$ and
%$\lim\limits_{n\to\be}f_n(x,y)=f(x,y)$ for all $(x,y)\in\mX\times\mY$.
%Thus, $f$ belongs to the first Baire class, and so it is Borel.
%>>>>>>>>>>>>>>>>>>>>
%
%\begin{thm}
%\label{36}
%Let $\mX$ and $\mY$ be closed subsets of $\C$ and let $\Phi$ be a
%function on $\mX\times\mY$ that is continuous in each variables.
%Suppose that $\mu$ and $\nu$ are positive regular Borel measures on $\mX$ and $\mY$
%such that $\supp\mu=\mX$ and $\supp\nu=\mY$.
%Then $\|\Phi\|_{\frak M_{\mX,\mY}}=\|\Phi\|_{\frak M_{\mX,\mY}^{\,\mu,\nu}}$.
%\end{thm}
The following result is also well known.

\medskip

{\it Let $f\in C(\C)$. Put $\Phi(z,w)\df f(z-w)$. Then $\Phi\in\frak
M_{\C}$ if and only if there exists a complex measure $\mu$ on $\C$ such that}
\bay
\label{kuda}
f(z)=\int_\C e^{-{\rm i}\re(z\ov\z)}\,d\mu(\z)
\quad\mbox{and}\quad\|\Phi\|_{\frak M_{\C}}=|\mu|(\C).
\ey

\medskip

A similar statement holds for an arbitrary locally compact abelian
group. The case of the group $\Z$ is considered, e.g., in \cite{Be}.

Let us proceed now to double operator integrals. The theory of
double operator integrals was developed by Birman and Solomyak in
\cite{BS1}, \cite{BS2}, and \cite{BS3}, see also their survey
\cite{BS4}.

Let $(\X,E_1)$ and $(\Y,E_2)$ be spaces with spectral measures $E_1$
and $E_2$ on a separable Hilbert space $\h$. The approach of Birman and
Solomyak is to define first double operator integrals
\bay
\label{doi}
\int\limits_\X\int\limits_\Y\Phi(x,y)\,d
E_1(x)T\,dE_2(y),
\ey
for bounded measurable functions $\Phi$ and
operators $T$ of Hilbert Schmidt class $\bS_2$.

We define here double operator integrals for arbitrary bounded operators $T$ and
refer the reader to \cite{BS1}, \cite{BS3}, and \cite{Pe2}.

\medskip

{\bf Definition.} Let $\mu$ and $\nu$ be $\s$-finite measures on $\mX$ and $\mY$
such that $E_1$ and $\mu$ are mutually absolutely continuous,
and $E_2$ and $\nu$ are mutually
absolutely continuous. We say that a measurable function $\Phi$ on $\mX\times\mY$
is {\it a Schur multiplier with respect to $E_1$ and $E_2$}
if $\Phi\in\fM_{\mX,\mY}^{\,\mu,\nu}$.
We denote the space of such Schur multipliers by $\fM(E_2,E_1)$.

\medskip

It is well known that the definition does not depend on the choice of
measures $\mu$ and $\nu$.

Let us now define double operator integrals \rf{doi} for bounded operators $T$.
Suppose that $\Phi\in\fM(E_2,E_1)$ and $\Phi$ admits a representation \rf{ipt} with
$$
\int_\O\|\f(\cdot,w)\|_{L^\be(E_1)}\|\psi(\cdot,w)\|_{L^\be(E_2)}\,d\l(w)<\be.
$$
We put
$$
\int\limits_\X\int\limits_\Y\Phi(x,y)\,dE_1(x)T\,dE_2(y)\df
\int\limits_\O\left(\,\int\limits_\X\f(x,w)\,dE_1(x)\right)T
\left(\,\int\limits_\Y\psi(y,w)\,dE_2(y)\right)\,d\l(w).
$$
It can be shown that the definition does not depend on the choice of
a representation \rf{ipt}.

It is also well known that for $\Phi\in\fM(E_2,E_1)$,
$$
\|\Phi\|_{\fM_{\mX,\mY}^{\,\mu,\nu}}=\|\Phi\|_{\fM(E_2,E_1)},
$$
where
$$
\|\Phi\|_{\fM(E_2,E_1)}\df\sup_{\|T\|\le1}\left\|
\int\limits_\X\int\limits_\Y\Phi(x,y)\,dE_1(x)T\,dE_2(y)\right\|
$$
and $\mu$ and $\nu$ are as in the above definition.

Birman and Solomyak proved that if $R$ is a bounded linear operator, $A$ and $B$ are
(not necessarily bounded) self-adjoint operators such that $AR-RB$ is bounded
and if $f$ is a continuously differentiable
function on $\R$ such that the divided difference $\dg f$ defined by
$$
(\dg f)(x,y)=\frac{f(x)-f(y)}{x-y}
$$
is a Schur multiplier with respect to the spectral measures of $A$ and $B$, then
\bay
\label{BSss}
f(A)R-Rf(B)=\iint_{\R\times\R}(\dg f)(x,y)\,dE_{A}(x)(AR-RB)\,dE_B(y)
\ey
and
$$
\|f(A)R-Rf(B)\|\le\|f\|_{\fM(E_A,E_{B})}\|AR-RB\|,
$$
(see \cite{BS3}).

Let us proceed now the case of normal operators. Suppose that
$N_1$ and $N_2$ are normal operators, $R$ is a bounded operator
such that the operator $N_1R-RN_2$ is bounded, $f$ is a continuous function on $\C$, and
the function $\dg_0f$ is defined by
$$
(\dg_0f)(\z_1,\z_2)\df\left\{\begin{array}{ll}
\frac{f(\z_1)-f(\z_2)}{\z_1-\z_2},&\z_1\ne\z_2,\\[.2cm]
0,&\z_1=\z_2.\end{array}\right.
$$
As in the case of self-adjoint operators, it can be shown that
\bay
\label{BSno}
f(N_1)R-Rf(N_2)=
\iint\limits_{\C\times\C}(\dg_0f)(\z_1,\z_2)\,dE_{N_1}(\z_1)(N_1R-RN_2)\,dE_{N_2}(\z_2),
\ey
whenever $\dg_0f\in\fM(E_{N_1},E_{N_2})$. Moreover,
$$
\|f(N_1)R-Rf(N_2)\|\le\|\dg_0f\|_{\fM(E_{N_1},E_{N_2})}.
$$
However, the class of functions $f$, for which formula \rf{BSno} can be used is not as ample in
general as in the case of formula \rf{BSss}. Indeed, it follows from results of \cite{JW} that if $N_1=N_2$
and the spectrum $\s(N_1)$ of $N_1$ has a nonisolated point, then $f$ must have complex
derivative at that point. In particular, $\dg_0f\in\fM(E_{N_1},E_{N_2})$
for all such normal operators $N_1$ and $N_2$ if and only if $f$ is a linear function.

\section{\bf Operator Lipschitz and commutator Lipschitz functions}
\setcounter{equation}{0}
\label{OLCL}

\

In this section we study properties of operator Lipschitz functions.
We also introduce the class of commutator Lipschitz functions.

We deal in this section  with bounded normal operators.
We show later that almost all the results remain
valid for unbounded normal operators, see \S\,\ref{Omc}.

Let $\SA$ denote the set of
bounded self-adjoint operators,  $\Un$ denote the set of unitary operators,
and let $\Pro$ denote the set of orthogonal projections.
For a closed subset $\fF$ of $\C$, we denote by $\No(\fF)$
the set of bounded normal operators $N$ with spectrum $\s(N)$ contained
in $\fF$;  $\No\df\No(\C)$. Finally, we put $\USA\df\Un\cap\SA$.

The following theorem is a generalization of Theorem 10.1 in \cite{AP2}.

\begin{thm}
\label{sr}
Let $f$ be a continuous function on a closed subset $\fF$ of $\C$. The following are equivalent:

{\em(i)} $\|f(N_1)-f(N_2)\|\le\|N_1-N_2\|$ for all $N_1, N_2\in\No(\fF)$;

{\em(ii)} $\|f(N)U-Uf(N)\|\le\|NU-UN\|$ for all $N\in\No(\fF)$ and
$U\in\Un$;

{\em(iii)} $\|f(N)A-Af(N)\|\le\|NA-AN\|$ for all $N\in\No(\fF)$
and $A\in \SA$;

{\em(iv)} $\|f(N)Q-Qf(N)\|\le\|NQ-QN\|$ for all $N\in\No(\fF)$ and $Q\in\USA$;

{\em(v)} $\|f(N)P-Pf(N)\|\le\|NP-PN\|$ for all $N\in\No(\fF)$ and $P\in\Pro$;

{\em(vi)} $\|f(N)R-Rf(N)\|\le\max\big\{\|NR-RN\|,\|N^*R-RN^*\|\big\}$ for all $N\in\No(\fF)$
and all bounded operators $R$;

{\em(vii)} $\|f(N_1)R-Rf(N_2)\|\le\max\big\{\|N_1R-RN_2\|,\|N_1^*R-RN_2^*\|\big\}$ for all $N_1, N_2\in\No(\fF)$
and all bounded operators $R$.
\end{thm}

\Pf The implications (vii)$\imp$(i) and (iii)$\imp$(iv) are trivial.
Note that $Q\in\USA$ if and only if $Q=2P-I$  for an orthogonal projection $P$.
Hence, statements (iv) and (v) are equivalent.

Thus it suffices to  verify the implications (i)$\imp$(ii)$\imp$(iii)$\imp$(vi)$\imp$(vii) and (iv)$\imp$(i).

To prove the implication (i)$\imp$(ii), it suffices to put $N_1\df N$ and $N_2\df UNU^*$.
Let us show that (ii)$\imp$(iii). Substituting $U=\exp({\rm i}tA)$ in (ii), we obtain
$$
\big\|f(N)-\exp({\rm i}tA)f(N)\exp(-{\rm i}tA)\big\|\le
\big\|N-\exp({\rm i}tA)N\exp(-{\rm i}tA)\big\|\quad\mbox{for every}\quad t\in\R.
$$
It remains to observe that
$$
\lim_{t\to0}\frac{\|f(N)-\exp({\rm i}tA)f(N)\exp(-{\rm i}tA)\|}{|t|}=\|f(N)A-Af(N)\|
$$
and
$$
\lim_{t\to0}\frac{\|N-\exp({\rm i}tA)N\exp(-{\rm i}tA)\|}{|t|}=\|NA-AN\|.
$$
Let us prove now that (iii)$\imp$(vi).
We consider the normal operator $\N$ and the bounded self-adjoint operator $\A$ defined
as follows
$$
\N=\left(\begin{matrix}N&\0\\[.2cm]\0&N\end{matrix}\right)\quad\mbox{and}\quad
\A=\left(\begin{matrix}\0&R\\[.2cm]R^*&\0\end{matrix}\right).
$$
It is easy to see that $\s(\N)=\s(N)$,
$$
f(\N)\A=\left(\begin{matrix}\0&f(N)R\\[.2cm]f(N)R^*&\0\end{matrix}\right)\quad\mbox{and}\quad
\A f(\N)=\left(\begin{matrix}\0&R f(N)\\[.2cm]R^*f(N)&\0\end{matrix}\right).
$$
Clearly,
$$
\|f(\N)\A-\A f(\N)\|=\max\big\{\|f(N)R-Rf(N)\|,~\|f(N)R^*-R^*f(N)\|\big\}
$$
and
$$
\|\N\A-\A\N\|=\max\big\{\|NR-RN\|,~\|NR^*-R^*N\|\big\}.
$$
It follows that
\begin{align*}
\|f(N)R-Rf(N)\|&\le\|f(\N)\A-\A f(\N)\|\le\|\N\A-\A\N\|\\[.2cm]
&=\max\big\{\|NR-RN\|,~\|NR^*-R^*N\|\big\}.
\end{align*}
Now let us show that (vi)$\imp$(vii). Put
$$
\N\df\left(\begin{matrix}N_1&\0\\[.2cm]\0&N_2\end{matrix}\right)\quad\mbox{and}\quad
\cR\df\left(\begin{matrix}\0&R\\[.2cm]\0&\0\end{matrix}\right).
$$
Then $\s(\N)=\s(N_1)\cup\s(N_2)$,
$$
f(\N)\cR=\left(\begin{matrix}\0&f(N_1)R\\[.2cm]\0&\0\end{matrix}\right),\quad\mbox{and}\quad
\cR f(\N)=\left(\begin{matrix}\0&Rf(N_2)\\[.2cm]\0&\0\end{matrix}\right).
$$
Hence, %$\|f(\N)\cR-\cR f(\N)\|=\|f(N_1)R-Rf(N_2)\|$.
\begin{align*}
\|f(N_1)R-Rf(N_2)\|&=\|f(\N)\cR-\cR f(\N)\|\\[.2cm]
&\le\max\big\{\|\N\cR-\cR\N\|,~\|\N\cR^*-\cR^*\N\|\big\}\\[.2cm]
&=\max\big\{\|NR-RN\|,~\|NR^*-R^*N\|\big\}.
\end{align*}

To complete the proof, it remains to show that (iv)$\imp$(i). Let $N_1,N_2\in\No(E)$. Put
$$
\N=\left(\begin{matrix}N_1&\0\\[.2cm]\0&N_2\end{matrix}\right)\quad\mbox{and}\quad
\mathscr Q=\left(\begin{matrix}\0&I\\[.2cm]I&\0\end{matrix}\right).
$$
Then $\N\in\No(E)$,
$$
f(\N)\mathscr Q=\left(\begin{matrix}\0&f(N_1)\\[.2cm]f(N_2)&\0\end{matrix}\right),\quad\mbox{and}\quad
\mathscr Q f(\N)=\left(\begin{matrix}\0&f(N_2)\\[.2cm]f(N_1)&\0\end{matrix}\right),
$$
and the inequality $\|f(\N)\mathscr Q-\mathscr Q f(\N)\|\le\|\N\mathscr Q-\mathscr Q \N\|$
implies the inequality \lb $\|f(N_1)-f(N_2)\|\le\|N_1-N_2\|$.
$\bl$

\medskip
The reasoning in the proof of (vi)$\imp$(vii) allows us to obtain
the following statement:

\begin{thm}
\label{sr0}
Let $f$ be a continuous function on a closed subset $\fF$ of $\C$. The following are equivalent:

{\em(i)} $\|f(N)R-Rf(N)\|\le\|NR-RN\|$ for all $N\in\No(\fF)$
and all bounded operators $R$;

{\em(ii)} $\|f(N_1)R-Rf(N_2)\|\le\|N_1R-RN_2\|$ for all $N_1, N_2\in\No(\fF)$
and all bounded operators $R$.
\end{thm}

Denote by $\OL(\fF)$ the set of operator Lipschitz functions on $\fF$, i.e., the set of $f\in C(\fF)$ such that
\bay
\label{lip}
\|f(N_1)-f(N_2)\|\le \const\|N_1-N_2\|\quad\mbox{for all}\quad N_1,\,N_2\in\No(\fF).
\ey
We use the notation $\|f\|_{\OL(\fF)}$ for the best constant on
the right-hand side of \rf{lip}.
It is easy to see that $\|f\|_{\OL(\fF)}\le1$ if and only if $f$ satisfies the equivalent
statements of Theorem \ref{sr}.

Denote by $\COL(\fF)$ the set of commutator Lipschitz functions on $\fF$, i.e., the set of $f\in C(\fF)$ such that
\bay
\label{comlip}
\|f(N)R-Rf(N)\|\le \const\|NR-RN\|%\quad\mbox{for all}\quad N\in\No(\fF).
\ey
for all $N\in\No(\fF)$ and bounded operators $R$.
We use the notation $\|f\|_{\COL(\fF)}$ for the best constant on
the right-hand side of \rf{comlip}. Theorem \ref{sr0} implies that
$$
\|f(N_1)R-Rf(N_2)\|\le \|f\|_{\COL(\fF)}\|N_1R-RN_2\|,
$$
whenever $N_1,\,N_2\in\No(\fF)$,  $R$ is a bounded operator, and
$f\in\COL(\fF)$.

It is clear that $\COL(\fF)\subset\OL(\fF)$ and $\|f\|_{\OL(\fF)}\le\|f\|_{\COL(\fF)}$.
Note also that $\COL(\fF)=\OL(\fF)$ if $\fF\subset\R$. But in general $\COL(\fF)\ne\OL(\fF)$.

For example, $\ov z\in\OL(\C)\setminus\COL(\C)$. Moreover, it is well known that
if $f\in\COL(\fF)$, then there exists finite limit
$$
\lim\limits_{z\to z_0}\dfrac{f(z)-f(z_0)}{z-z_0}
$$
for each limit point $z_0$ of $\fF$. This follows from results of
\cite{JW}, see also \cite{KS}.
Indeed, inequality \rf{ssylka} below implies that
$\dg_0f\in{\frak M}_{\frak F}$. Hence, $f$ has complex derivative
at any nonisolated point of $\frak F$ by Theorem 4.1 in \cite{JW}.

In particular, $\COL(\C)=\{az+b:a,b\in\C\}$.

\begin{thm}
\label{diflip}
Let $f$ be an operator Lipschitz function on $\C$.
Then $f$ has finite derivative at every point in every direction.
\end{thm}

\Pf Clearly, $f\big|\R$ is an operator Lipschitz function on $\R$.
Hence, it is differentiable everywhere on $\R$ by Theorem 4.1 in \cite{JW}.
In particular, the
partial derivative $\frac{\partial f}{\partial x}(0)$ exists and is
finite. To complete the proof, it suffices to observe that the space
of operator Lipschitz functions on $\C$ is translation and rotation
invariant. $\bl$

\medskip

It is clear from the proof of Theorem \ref{diflip} that the following statement
is also true.

\begin{thm}
Let $f\in\OL(\fF)$, where $\fF$ be a closed subset of $\C$.
Then $f$ has finite derivative at every point $\z\in\fF$ in each direction
$\xi\in\C$ such that $0$ is not an isolated point of $\{t\in\R:\z+t\xi\in\fF\}$.
\end{thm}

Nevertheless, it turns out that functions in $\OL(\C)$ do not have to be differentiable
as functions of two real variables.
To construct such a function, we put
$$
h_n(\z)\df\left\{\begin{array}{ll}\z^{n+1}/\bar\z^n,&\text {if}\,\,\,\,\z\not=0,\\[.2cm]
0,&\text {if}\,\,\,\,\z=0.
\end{array}\right.
$$

\begin{thm}
\label{hn}
Let $n\in\Z$.
Then $h_n\in\OL(\C)$ and
$$
\|h_n\|_{\OL(\C)}=\|h_n\|_{\Li(\C)}=|2n+1|.
$$
The function $h_n$ is not differentiable at the origin unless $n=0$ or $n=-1$.
\end{thm}

\Pf Clearly, $h_0(\z)=\z$ and $h_{-1}(\z)=\bar\z$. Thus the conclusion of the theorem is obvious if $n=0$ or $n=-1$.
Put
$$
\sgn\z\df\left\{\begin{array}{ll}\z/|\z|,&\text {if}\,\,\,\,\z\not=0,\\[.2cm]
0,&\text {if}\,\,\,\,\z=0.
\end{array}\right.
$$
It is easy to see that $h_n(\z)\df\z\sgn^{2n}\z$ for $n>0$ and
 $\ov h_n=h_{-n-1}$ for every $n$ in $\Z$. Hence, it suffices to consider the case $n>0$.

Let $N_1$ and $N_2$ be normal operators. We have
\begin{align*}
h_n(N_1)-h_n(N_2)=&\sum_{j=0}^n\sgn^{2n-2j}(N_1)(N_1-N_2)\sgn^{2j}(N_2)\\[.2cm]
&+\sum_{j=0}^{n-1}\sgn^{2n-2j}(N_1)(N_2^*-N_1^*)\sgn^{2j+2}(N_2),
\end{align*}
since
$$
N_2\sgn^{2j}(N_2)=N_2^*\sgn^{2j+2}(N_2)
$$
and
$$
\sgn^{2n-2j-2}(N_1)N_1=\sgn^{2n-2j}(N_1)N_1^*.
$$
Hence, $\|h_n(N_1)-h_n(N_2)\|\le (2n+1)\|N_1-N_2\|$.

Since, obviously, $\|h_n\|_{\OL(\C)}\ge\|h_n\|_{\Li(\C)}$, it suffices to show that $\|h_n\|_{\Li(\C)}\ge2n+1$, which follows immediately from the equality
$h_n(e^{{\rm i}t})=e^{{\rm i}(2n+1)t}$.

It is easy to see that $h_n$ is not differentiable at the origin
for $n\notin\{0,-1\}$. $\bl$

For a function $f$ on a subset $\fF$ of $\C$ we define the function
$$
\big(\dg_0f\big)(z,w)\df
\left\{\begin{array}{ll}\frac{f(z)-f(w)}{z-w},&\text {if}\,\,\,\,z,w\in \frak F, \,\,\,z\ne w,\\[.2cm]
0,&\text {if}\,\,\,\,x\in \frak F, \,\,\,z=w.
\end{array}\right.
$$
We need the following well-known inequality \bay \label{nva}
\|f\|_{{\rm CL}(\frak F)}\le\|\dg_0 f\|_{{\frak M}_{\frak F}} \ey
for any closed subset $\fF$ of $\C$. In the case $\fF\subset\R$ the proof
can be found in \cite{AP6}. Inequality \rf{nva} can be derived from the
following formula (see \rf{BSno}):
$$
f(N_1)R-Rf(N_2)=\iint\limits_{\fF\times\fF}(\dg_0
f\big)(z,w)\,dE_{N_1}(z)(N_1R-RN_2)\,dE_{N_2}(z),
$$
where $f$ is a function such that $\dg_0f\in{\frak M}_{\frak F}$,
$N_1$ and $N_2$ are normal operators with bounded $N_1R-RN_2$ whose
spectra are in $\fF$, and $E_{N_1}$ and $E_{N_2}$ are the spectral
measures of $N_1$ and $N_2$. One can prove that
\bay
\label{ssylka}
\|\dg_0 f\|_{{\frak M}_{\frak F}}\le2\|f\|_{{\rm CL}(\frak F)}
\ey
for every closed subset $\fF$ of $\C$. This was proved in \cite{AP6}
in the case $\fF\subset\R$. The general case can be treated in the
same way. We omit the details because we are not going to apply this
estimate in this paper.

Let $f\in{\rm CL}(\frak F)$. Suppose that $\frak F$ has no isolated
points. Put
$$
\big(\dg f\big)(z,w)\df
\left\{\begin{array}{ll}\frac{f(z)-f(w)}{z-w},&\text {if}\,\,\,\,x,¸áy\in \frak F, \,\,\,z\ne w,\\[.2cm]
\lim\limits_{\z\to w}\frac{f(\z)-f(w)}{\z-w},&\text {if}\,\,\,\,z\in \frak F, \,\,\,z=w.
\end{array}\right.
$$
It was observed in \cite{AP6} that the equality
$$
\|f\|_{{\OL}(\fF)}=\|f\|_{{\rm CL}(\fF)}=\|\dg f\|_{{\fM}_{\fF}}
$$
holds in the case $\fF\subset\R$. In the same way one can prove that
$$
\|f\|_{{\COL}(\fF)}=\|\dg f\|_{{\fM}_{\fF}}
$$
for every closed subset $\fF$ of $\C$. We omit the details.

As we have mentioned above, the operator Lipschitz norm admits the following
estimate in terms of the multiplier norm:
\bay
\label{D0f}
\|f\|_{\OL(\fF)}\le\|\dg_0 f\|_{{\frak M}_{\frak F}}.
\ey
However, this estimate is not as helpful as in the commutator Lipschitz case.
Indeed, if $\fF$ has nonempty interior, then for the function $\ov z$
the right-hand side of \rf{D0f} is infinite, while the function $\ov z$ is
certainly operator Lipschitz.

On the other hand, the operator Lipschitz norm admits the following upper
estimate in terms of the multiplier norms of certain functions:

{\it if $f$ is a continuous function on a closed subset $\fF$ of $\C$ that
admits a representation
\bay
\label{kuzya}
f(z)-f(w)=(z-w)g_1(z,w)+(\ov z-\ov w)g_2(z,w),\quad z,w\in\fF,
\ey
for some $g_1,g_2\in\fM_{\fF}$. Then $f\in\OL(\fF)$ and}
\bay
\label{dunya}
\|f\|_{\OL(\fF)}\le\|g_1\|_{\fM(\fF)}+\|g_2\|_{\fM(\fF)}.
\ey

Indeed, as in Theorems 5.2 and 10.1 of \cite{APPS2}, it can be shown that
the following formula holds for $N_1,\,N_2\in\No(\fF)$:
\begin{align*}
f(N_1)-f(N_2)&=\iint\limits_{\fF\times\fF} g_1(z,w)\,dE_{N_1}(z)(N_1-N_2)\,dE_{N_2}(w)\\[.2cm]
&+\iint\limits_{\fF\times\fF} g_2(z,w)\,dE_{N_1}(z)(N^*_1-N^*_2)\,dE_{N_2}(w).
\end{align*}
This immediately implies \rf{dunya}.

Note that estimate \rf{dunya} in many cases is considerably more helpful than \rf{D0f}.
In particular, it was shown in \cite{APPS2} that if $f$ belongs to the Besov class
$B_{\be1}^1(\R^2)$, then $f$ admits a representation of the form \rf{kuzya} with
$g_1,\,g_2\in\fM(\C)$.

\

\section{\bf Absolutely convergent Fourier integrals\\
and estimates of commutator Lipschitz norms}
\setcounter{equation}{0}
\label{FouIn}

\

We are going to obtain in this section norm estimates of certain functions in
the space of functions with absolute convergent Fourier integrals. This will allow us
to obtain certain commutator Lipschitz estimates.

We denote by $\mZ$ the set of complex integers:
$$
\mZ\df\Z+{\rm i}\Z.
$$

Put
$$
\widehat L^1=\widehat L^1(\C)\df\mathscr F(L^1(\C)),\quad
\|f\|_{\widehat L^1}=\|f\|_{\widehat L^1(\C)}\df\big\|\mathscr F^{-1} f\big\|_{L^1}.
$$

\begin{lem}
\label{log0}
Suppose that $0<\d<{r}<\be$. Then there exists a function $h\in \widehat L^1(\C)$
such that
$$
\|h\|_{\widehat L^1}\le\const\log\frac{2{r}}\d,\quad h(0)=0,\quad
\mbox{and}\quad h(z)=\frac{\ov z}{z},\quad\mbox{whenever}\quad
\d\le|z|\le{r}.
$$
\end{lem}

\Pf Clearly, it suffices to consider the case when $\d=2$ and ${r}=2^n$ with $n\ge1$.
We can take an even function $\f\in C^\be(\R)$ such that $\supp\f\subset[-2,2]$
and $\f(x)=1$ for $x\in[-1,1]$. Put
$h(z)\df\dfrac{\ov z}{z}\big(\f(2^{-n}|z|)-\f(|z|)\big)$. Then $h\in\widehat L^1(\C)$
and $h(z)=\dfrac{\ov z}{z}$ in the annulus $\{2\le|z|\le2^n\}$.
To estimate $\|h\|_{\widehat L^1}$, we put $\psi(z)\df\dfrac{\ov z}{z}\big(\f(|z|/2)-\f(|z|)\big)$. Then
$$
h(z)=\sum_{k=0}^{n-1}\psi(2^{-k}z)
$$
and
$$
\|h\|_{\widehat L^1}\le\sum_{k=0}^{n-1}\big\|\psi(2^{-k}z)\big\|_{\widehat L^1}=n\|\psi\|_{\widehat L^1}.\quad\bl
$$

\begin{cor}
\label{log2}
Let $\L$ be a finite subset of $\C$. Assume that
$0<\d\le|\l-\mu|\le{r}$ for all $\l,\mu\in\L$ such that $\l\ne\mu$. Then
for $f(z)=\ov z$ we have
$$
\|f\|_{\COL(\L)}\le\|\dg_0 f\|_{\fM_{\L}}\le\const\log\frac{2{r}}\d.
$$
\end{cor}

\Pf The left inequality is a special case of \rf{nva} while the right
one is an immediate consequence of \rf{kuda} and Lemma \ref{log0}. $\bl$

\begin{cor}
\label{log1}
Let $f(z)=\ov z$ and $0<\d<{r}$. Then
$$
\|f\|_{\COL(\d\mZ\cap\,{r}\dd)}\le\|\dg_0 f\|_{\fM_{\d\mZ\cap\,{r}\dd}}
\le\const\log\frac{2{r}}\d,
$$
where $\dd\df\{z\in\C:|z|<1\}$.
\end{cor}

\begin{lem}
\label{ovz}
Put
\bay
\label{Phi}
\Psi(z)\df\left\{\begin{array}{ll}\ov z,&\text {if}\,\,\,|z|<1,\\[.2cm]
z^{-1},&\text {if}\,\,\,\,|z|\ge1.
\end{array}\right.
\ey
Then $\Psi\in \widehat L^1(\C)$.
\end{lem}

\Pf It is easy to see that $\frac{\partial\Psi}{\partial\ov z}=\chi_{_{\dd}}$ in the sense
of distributions. We need the well-known formula
\bay
\label{J1}
(\F^{-1}\chi_{_{\dd}})(\z)=\frac{1}{2\pi |\z|}J_{1}(|\z|),
\ey
where $J_1$ denotes the Bessel function. We prove \rf{J1} here for the reader's convenience.
Recall that
$$
J_\nu(x)\df\sum\limits_{k=0}^\infty\frac{(-1)^k(x/2)^{2k+\nu}}{\G(\nu+k+1)k!}, \quad\nu\in\C.
$$
Applying the polar change of variables and the Poisson formula (see \cite{W}, \S 2.3,
formula (2)) we find that
\begin{align*}
(\F^{-1}\chi_{_{\dd}})(\z)&
\df\frac1{4\pi^2}\int_\C\chi_{_{\dd}}(\xi)\,e^{\rm i\re(\xi\ov\z)}\,d\m_2(\xi)\\[.2cm]
&=
\frac1{4\pi^2}\int_0^1 r\int_{-\pi}^\pi e^{\rm i r|\z|\cos(\psi-\theta)}\,d\psi\\[.2cm]
&=\frac1{2\pi}\int_0^1 rJ_0(r|\z|)\,dt=\frac{1}{2\pi |\z|}J_{1}(|\z|),
\end{align*}
where $\m_2$ stands for planar Lebesgue measure.
Hence,
$$
(\F^{-1}\Psi)(\z)=\frac{2{\rm i}(\F^{-1}\chi_{_{\dd}})(\z)}{\z}
=\frac{{\rm i}J_1(|\z|)}{\pi\z|\z|}.
$$
It remains to observe that $|J_1(x)|\le\const x^{-1/2}$, $x>1$,
see, for example, \cite{W}, \S 7.21. $\bl$

\begin{cor}
\label{corovz}
The function
$$
(z,w)\mapsto\Psi\left(\frac{z-w}{a}\right)
$$
belongs to $\fM_\C$ and its norm in $\fM_\C$ is equal to
$\|\Psi\|_{\widehat L^1(\C)}$ for every $a>0$.
\end{cor}

Consider the following function on $\C$:
$$
\l(\z)\df
\left\{\begin{array}{ll}\z^{-1},&\z\ne 0,\\[.1cm]
0,&\z=0.
\end{array}\right.
$$
Note that for every function $f$ on a subset $\fF$ of $\C$,
$$
(\dg_0 f)(z,w)=\big(f(z)-f(w)\big)\l(z-w),\quad z,\,w\in\fF.
$$

\begin{thm}
Let $\L$ be a subset of $\C$. Suppose that $|\l-\mu|\ge\d>0$
for all distinct $\l$ and $\mu$ in $\L$. Then
$$
\|f\|_{\COL(\L)}\le\|\dg_0 f\|_{\fM_{\L}}\le\const \frac1\d\sup\big\{|f(\l)|:\l\in\L\big\}.
$$
\end{thm}

\Pf
Note that
$$
\l(z-w)=\frac1\d\Psi\left(\frac{z-w}{\d}\right),\quad(z,w)\in \L\times \L,
$$
where $\Psi$ is defined by \rf{Phi}.
Hence,  by Corollary \ref{corovz},
$$
\|\l(z-w)\|_{\fM_{\L}}\le\const\frac1\d,
$$
and so
\begin{align*}
\|\dg_0 f\|_{\fM_{\L}}&\le\|f(z)\l(z-w)\|_{\fM_{\L}}+
\|f(w)\l(z-w)\|_{\fM_\L}\\[.2cm]
&\le2\const \frac1\d\sup\big\{|f(\l)|:\l\in\L\big\}.\quad\bl
\end{align*}

\begin{cor}
\label{lat}
Let $\d>0$ and let $f$ be a bounded function on $\d\mZ$. Then
$$
\|f\|_{\COL(\d\mZ)}\le\|\dg_0 f\|_{\fM_{\,\d\mZ}}
\le\const \frac1\d\sup\big\{|f(\d\fn)|:\fn\in\C\big\}.
$$
\end{cor}

The following theorem shows that Corollary \ref{log1} is sharp.

\begin{thm}
\label{niz}
Let $f(z)=\ov z$ and $0<\d<{r}$. Then
$$
\|f\|_{\COL(\d\mZ\cap\,{r}\dd)}\ge\const\log\frac{2{r}}\d.
$$
\end{thm}

To prove the theorem, we need several auxiliary facts.

\begin{lem}
\label{52}
Let $\Psi$ be the function defined by {\em\rf{Phi}}. Then
$$
|((\mathscr F^{-1}(\Psi^2))(\xi)|\le\frac{\const}{1+|\xi|^{\frac52}},\quad\xi\in\C.
$$
\end{lem}

\Pf It is easy to see that $\frac{\partial\Psi^2}{\partial\ov z}=2\ov z\chi_{_{\dd}}(z)$ in the sense
of distributions.

Applying the polar change of variables, we obtain
\begin{align*}
(\F^{-1}(2\ov \xi\chi_{_{\dd}}(\xi)))(\z)&\df
\frac1{2\pi^2}\int_\C\ov\xi\,\chi_{_{\dd}}(\xi)\,e^{\rm i\re(\xi\ov\z)}\,d\m_2(\xi)\\[.2cm]
&=
\frac1{2\pi^2}e^{-\rm i\theta}\int_0^1 r^2
\int_{-\pi}^\pi e^{\rm i r|\z|\cos(\psi-\theta)}e^{-\rm i(\psi-\theta)}\,d\psi\\[.2cm]
&=\frac1{\pi^2}e^{-\rm i\theta}\int_0^1 r^2\int_{0}^\pi e^{\rm i r|\z|\cos\psi}\cos\psi\,d\psi,
%=\frac1{2\pi}\int_0^1 rJ_0(r|\z|)\,dt=\frac{1}{2\pi |\z|}J_{1}(|\z|),
\end{align*}
where $\m_2$ stands for planar Lebesgue measure.
Applying the Poisson formula (see \cite{W}, \S 2.3,
formula (2)), we obtain
$$
\int_0^\pi e^{\rm i x\cos\psi}\cos\psi\,d\psi=-{\rm i}\frac d{dx}\left(\int_0^\pi e^{\rm i x\cos\psi}\,d\psi\right)
=-\rm i \pi J_0^\prime(x)=\rm i \pi J_1(x).
$$
Hence,
\begin{align*}
(\F^{-1}(2\ov \xi\chi_{_{\dd}}(\xi)))(\z)&=
\frac{{\rm i}e^{-\rm i\theta}}\pi \int_0^1 r^2 J_1(r|\z|)\,dr
=\frac{{\rm i} e^{-\rm i\theta}}{\pi|\z|^3}\int_0^{|\z|}t^2J_1(t)\,dt\\[.2cm]
&=\frac{{\rm i} e^{-\rm i\theta}}{\pi|\z|^3}\int_0^{|\z|}(t^2J_2(t))^\prime\,dt
=\frac{{\rm i} e^{-\rm i\theta}J_2(|\z|)}{\pi|\z|},
\end{align*}
and so
$$
(\mathscr F^{-1}(\Psi^2))(\z)=-\frac{2 e^{-\rm i\theta}J_2(|\z|)}{\pi\z|\z|}=-\frac{2 J_2(|\z|)}{\pi\z^2}.
$$
It remains to observe that $|J_2(x)|\le\const x^2$ for $x\in[0,1]$ and $|J_2(x)|\le\const x^{-1/2}$
for $x\in[1,\be)$;
see, for example, \cite{W}, \S 7.21. $\bl$

\begin{cor}
\label{blin}
The matrix $\{\l^2(\fm-\fn)\}_{\fm,\fn\in\mZ}$ induces
a bounded operator on $\ell^2(\mZ)$.
\end{cor}

\Pf Put
$$
h(\z)=\sum_{\fn\in\mZ\setminus\{0\}}\frac{e^{\rm i\re(\z\ov\fn)}}{\fn^2}, \quad\z\in\C.
$$
Clearly, the series converges in $L^2([0,2\pi]^2)$ and $h(\z)=h(\z+2\pi)=h(\z+2\pi\rm i)$.
It suffices to verify that $h\in L^\be(\C)$. Put
$$
h_0(\z)\df4\pi^2\sum_{\fn\in\Z}H(\z+2\pi\fn),
$$
where $H\df\F^{-1}(\Psi^2)$. Lemma \ref{52} implies that $h_0\in L^\be(\C)$. Moreover,
$$
\widehat h_0(\fn)\df\frac1{4\pi^2}\int_{[0,2\pi]^2}h_0(\xi)e^{-\rm i\re(\xi\ov\fn)}\,dm_2(\xi)
=(\F H)(\fn)=\Psi^2(\fn)
%=\frac1{\fn^2}
=\l^2(\fn)=\widehat h(\fn)
$$
for all $\fn\in\mZ$. Hence, $h=h_0$ almost everywhere on $\C$. $\bl$

For $r>0$, consider the matrix
\bay
\label{Lar}
\L_{r}\df\{|\l(\fm-\fn)|^2\}_{\fm,\fn\in\mZ\cap\,{r}\dd}.
\ey
The following lemma gives a lower estimate for the operator norm $\|\L_{r}\|$ of
$\L_r$.

\begin{lem}
\label{rho}
Let ${r}\ge3$. Then $\|\L_{r}\|\ge\const\log{r}$.
\end{lem}

\Pf Let us first observe that if $A=\{a_{jk}\}_{1\le j,k\le n}$
is a matrix with nonnegative
entries and $v=\{v_j\}_{1\le j\le n}$ is the vector with
$v_j=1$, $1\le j\le n$, then
$$
\sum_{j,k=1}^na_{jk}=(Av,v)\le n\|A\|.
$$
Thus it suffices to prove that
$$
\sum_{\fm,\fn\in\mZ\cap\,{r}\dd}|\l(\fm-\fn)|^2\ge\const{r}^2\log{r}.
$$
Put
$$
a_\fp\df\card\{(\fm,\fn)\in\mZ^2:\fm-\fn=\fp, |\fm|\le{r}, |\fn|\le{r}\},
$$
where $\fp\in\mZ$.
Clearly, $a_\fp\ge\const{r}^2$ if $|\fp|\le\frac12{r}$. We have
$$
\sum_{\fm,\fn\in\mZ\cap{r}\ov\dd}|\l(\fm-\fn)|^2=\sum_{\fp\in\mZ,\,\,\, \fp\ne0}\frac{a_\fp}{|\fp|^2}
\ge\const{r}^2\sum_{\fp\in\mZ,\,\,\,0<2|\fp|\le{r}}\frac{1}{|\fp|^2}
\ge\const{r}^2\log{r}.\quad \bl
$$

\medskip

{\bf Proof of Theorem \ref{niz}.} The theorem can be reformulated as follows:
$$
\|\dg_0f\|_{\fM_{\d\mZ\cap\,{r}\dd}}\ge\const\log\frac{2{r}}\d.
$$
Clearly, it suffices to assume that $\d=1$ and ${r}\ge3$.
Consider the matrices
$$
\L^{[1]}_{r}\df\{\dg_0f(\fm,\fn)\}_{\fm,\fn\in\mZ\cap\,{r}\dd}\quad\mbox{and}\quad
\L^{[2]}_{r}\df\{\ov\l\,^2(\fm-\fn)\}_{\fm,\fn\in\mZ\cap\,{r}\dd}.
$$
We have
$$
\L^{[1]}_{r}\star\L^{[2]}_{r}=\L_{r},
$$
where the matrix $\L_{r}$ is defined by \rf{Lar}
and the Schur--Hadamard product of matrices
is defined by \rf{ScHad}.

It remains to observe that $\big\|\L^{[2]}_{r}\big\|\le\const$
by Corollary \ref{blin} and $\|\L_{r}\|\ge\const\log{r}$
by Lemma \ref{rho}. $\bl$

\

\section{\bf Operator and commutator moduli of continuity}
\setcounter{equation}{0}
\label{Omc}

\

In this section we define the operator modulus of continuity and various
versions of commutator modulus of continuity. We study their properties and
obtain estimates that will be used in the next section.

Let $f$ be a continuous function defined on a closed subset $\fF$ of $\C$. Put
$$
\O_{f}(\d)=\O_{f,\fF}(\d)\df\sup\big\|f(N_1)-f(N_2)\big\|,\quad\d>0,
$$
where the supremum is taken over all $N_1,N_2\in\No(\fF)$ such that $\|N_1-N_2\|\le\d$. We say that $\O_f$ is the
{\it operator modulus of continuity of} $f$.

If $f$ is defined on a set that contains $\fF$, we put
$\O_{f,\fF}\df\O_{f|\fF,\fF}$.

The case $\fF=\R$ was considered in \cite{AP2}.

Note that the function $f$ on $\fF$ is operator Lipschitz if and only if $\O_f(\d)\le\const\d$, $\d>0$.

Clearly, for every $f\in C(\fF)$,
$$
\o_f(\d)\le\O_f(\d)\,\quad\d>0,
$$
where $\o_f$ is the (scalar) modulus of continuity of $f$:
$$
\o_f(\d)\df\sup\{|f(\z_1)-f(\z_2)|:~|\z_1-\z_2|\le\d\},\quad\d>0.
$$
On the other hand, it was proved in \cite{APPS2}, Theorem 8.2 that for $f\in C(\C)$,
$$
\O_f(\d)\le\const\,(\o_f)_*(\d),
$$
where for a modulus of continuity $\o$, the functions $\o_*$ is defined by \rf{o*}.
In fact, the same is true for an arbitrary closed subset $\fF$ of $\C$.

\begin{thm}
\label{W}
Let $\o$ be a modulus of continuity and
let $\fF$ be a closed subset of $\C$. Then for every
$f\in\L_\o(\fF)$,
$$
\O_{f,\fF}\le C\|f\|_{\L_\o(\fF)}\,\o_*(\d),
$$
where $C$ is a numerical constant.
\end{thm}

\Pf The result reduces to the case $\fF=\C$ because each function $f\in\L_\o(\fF)$
extends to a function $f_\natural\in\L_\o(\C)$
so that $\|f_\natural\|_{\L_\o(\C)}\le\const\|f\|_{\L_\o(\fF)}$.
The appropriate extension can be constructed with the help of Whitney
type theorems, see \cite{S} for details. $\bl$

Let $f$ be a continuous function on a closed subset $\fF$ of $\C$. For $\d>0$, put
$$
\O_{f}^{\SA}(\d)\df\sup\big\{\|f(N)A-Af(N)\|:~N\in\No(\fF),
~A\in\SA,~\|A\|=1,~\|NA-AN\|\le\d\big\}
$$
and
$$
\O_{f}^{\CO}(\d)\df\sup\big\{\|f(N)R-Rf(N)\|:~N\in\No(\fF),
~\|R\|=1,~\|NR-RN\|\le\d\big\}.
$$

As in the case of $\O_f$, we can write $\O_{f,\fF}^{\SA}$ and $\O_{f,\fF}^{\CO}$
if we want to emphasize the dependence on a closed set $\fF$.

Note that in \cite{AP2} and \cite{AP4} in the case of subsets of the real
line
the notation $\O_f^{\flat}$ (and $\O_{f,\fF}^{\flat}$) was used for
$\O_{f,\fF}^{\SA}$
and $\O_{f,\fF}^{\CO}$.

Recall that
$\O_{f,\fF}^{\SA}=\O_{f,\fF}^{\CO}$ if $\fF\subset\R$.
This was proved in \cite{AP2} in the case $\fF=\R$ and it was observed in \cite{AP4}
that the same reasoning works in the case $\fF\subset\R$.

However, the equality does not hold for arbitrary subsets of $\C$.
For example, if $f(z)=\ov z$, then
$\O_{f,\C}^{\SA}(\d)=\d$ and $\O_{f,\C}^{\CO}(\d)=\be$.
The first equality is trivial.
To prove the second equality, we observe that
$\O_{f,\C}^{\CO}(\d)=
\d\O_{f,\C}^{\CO}(1)$ because $f$ is a homogeneous function of
degree 1. Thus,
$\O_{f,\C}^{\CO}(\d)=\be$ if and only if $f$ is not commutator Lipschitz.
The fact that $f$ is not commutator Lipschitz
follows from Corollary 4.3 in \cite{JW}.

\medskip

{\bf Remark.} It is easy to see that $\O_f^{\SA}=\O_{\ov f}^{\SA}$ for every $f\in C(\C)$.
However, as we have already observed, $\O_f^{\CO}\ne\O_{\ov f}^{\CO}$ for $f(z)=z$.

\medskip

Recall that an operator $R$ on Hilbert space is called a {\it contraction}
if $\|R\|\le1$.

The following two theorems show that Theorem 10.2 in \cite{AP2}
can be generalized to the case of normal operators with spectrum in
a fixed closed subset $\fF$ of $\C$.

\begin{thm}
\label{123}
Let $f$ be a continuous function on a closed subset $\fF$ of $\C$. Then
$$
\|f(N_1)R-Rf(N_2)\|\le\O_{f,\fF}^{\SA}\big(\max\{\|N_1R-RN_2\|,\|N_1^*R-RN_2^*\|\}\big)
$$
for arbitrary $N_1, N_2\in\No(\fF)$ and contractions $R$.
\end{thm}

\Pf We consider first the case $N_1=N_2$. Replacing $R$ by $R+\a I$
with $\|R+\a I\|=1$, we see that the case $\|R\|<1$ is reduces to the case $\|R\|=1$.
Now the desired inequality in the case $N_1=N_2$ can be proved in the same way as the implication (iii)$\imp$(vi)
in the proof of Theorem \ref{sr}.
The general case can be reduced to the case $N_1=N_2$ as in the proof of
the implication
(vi)$\imp$(vii) in Theorem \ref{sr}. $\bl$

\medskip

{\bf Remark.} Theorem \ref{123} shows that the definition of $\O_f^{\SA}$ can be
replaced with the following one:
$$
\O_{f}^{\SA}(\d)=\sup\big\{\|f(N_1)A-Af(N_2)\|\},
$$
where the supremum is taken over $N_1$ and $N_2$ in $\No(\fF)$ and self-adjoint
operators $A$ with $\|A\|=1$ such that $\|N_1A-AN_2\|\le\d$.

\begin{thm}
\label{mcc}
Let $f$ be a continuous function on a closed subset $\fF$ of $\C$. Then
$$
\O_f\le\O_f^{\SA}\le2\O_f.
$$
\end{thm}

\Pf The inequality $\O_f\le\O_f^{\SA}$ follows from Theorem \ref{123}.
To prove the inequality $\O_f^{\SA}\le2\O_f$, we repeat the arguments of
the corresponding part of the proof Theorem 10.2 in \cite{AP2}.
In particular, we will use the following
inequality (see \cite{AP2}, Lemma 10.4)
\bay
\label{vl}
\big\|(I-T^2)^{1/2}X-X(I-T^2)^{1/2}\big\|\le\frac{\|T\|\cdot\|XT-TX\|}{(1-\|T\|^2)^{1/2}},
\ey
where $X$ is a bounded operator and $T$ is a self-adjoint operator with $\|T\|<1$.

Let $A$ be a self-adjoint operator with $\|A\|=1$ and $\t\in(0,1)$.
Consider the operators
$$
\N=\left(\begin{matrix}N&\0\\[.2cm]\0&N\end{matrix}\right)\quad\mbox{and}\quad
\cU=\left(\begin{matrix}\t A&(I-\t^2A^2)^{1/2}\\[.2cm]-(I-\t^2A^2)^{1/2}&\t A\end{matrix}\right).
$$
Clearly, $\cU$ is a unitary operator. We have
$$
f(\N)\cU=\left(\begin{matrix}\t f(N)A&f(N)(I-\t^2A^2)^{1/2}\\[.2cm]-f(N)(I-\t^2A^2)^{1/2}&\t f(N)A\end{matrix}\right)
$$
and
$$
\cU f(\N)=\left(\begin{matrix}\t Af(N)&(I-\t^2A^2)^{1/2}f(N)\\[.2cm]-(I-\t^2A^2)^{1/2}f(N)&\t Af(N)\end{matrix}\right).
$$
Obviously,
$$
\|f(\N)\cU-\cU f(\N)\|\ge\t\|f(N)A-Af(N)\|.
$$
Applying \rf{vl} with $X=N$ and $T=\t A$, we find that
\begin{align*}
\|\N\cU-\cU\N\|&\le\t\|NA-AN\|+\big\|N\big(I-\t^2A^2\big)^{1/2}-\big(I-\t^2A^2\big)^{1/2}N\big\|\\[.2cm]
&\le\big(\t+\t^2(1-\t^2)^{-1/2}\big)\|NA-AN\|.
\end{align*}
Hence,
\begin{align*}
\|f(N)A-Af(N)\|&\le\t^{-1}\|f(\N)\cU-\cU f(\N)\|=\t^{-1}\|f(\cU^*\N\cU)-f(\N)\|\\[.2cm]
&\le\t^{-1}\O_f\big(\big\|\cU^*\N\cU-\N\big\|\big)=\t^{-1}\O_f\big(\big\|\N\cU-\cU\N\big\|\big)\\[.2cm]
&\le\t^{-1}\O_f\Big(\big(\t+\t^2(1-\t^2)^{-1/2}\big)\|NA-AN\|\Big).
\end{align*}
Taking $\t=1/2$, we obtain
\begin{align*}
\|f(N)A-Af(N)\|&\le2\O_f\left(\left(\frac12+\frac1{2\sqrt{3}}\right)\|NA-AN\|\right)\le2\O_f\big(\|NA-AN\|\big).\quad\bl
\end{align*}

%Put $\O_f^{\SA}\df\O_f^{[1]}=\O_f^{[2]}=\O_f^{[3]}$.

\medskip

{\bf Remark.}
Note that in general for continuous functions $f$ on $\C$, $\O_f\ne\O_f^{\SA}$.
Indeed, it was shown in \cite{AP4} that there are continuous functions $f$ on $\R$
such that \lb$\O_{f,\R}\ne\O_{f,\R}^{\SA}$.
On the other hand, it is easy to see that if $f$ is a continuous function on
$\R$ and $F(\z)=f(\re\z)$, $\z\in\C$, then  $\O_{F,\C}=\O_{f,\R}$ and $\O_{F,\C}^{\SA}=\O_{f,\R}^{\SA}$.

\medskip

\begin{thm}
\label{com123}
Let $f$ be a continuous function on a closed subset $\fF$ of $\C$. Then
$$
\|f(N_1)R-Rf(N_2)\|\le\O_{f,\fF}^{\CO}(\|N_1R-RN_2\|)
$$
for arbitrary $N_1, N_2\in\No(\fF)$ and arbitrary contractions $R$.
\end{thm}

\Pf The proof is similar to the proof of Theorem \ref{123}. $\bl$

\medskip

{\bf Remark.}
Theorem \ref{com123} shows that the definition of $\O_f^{\CO}$ can be
replaced with the following one:
$$
\O_{f}^{\CO}(\d)=\sup\big\{\|f(N_1)R-Rf(N_2)\|:~N_1,\,N_2\in\No(\fF),
~\|R\|=1,~\|N_1R-RN_2\|\le\d\big\}.
$$

\begin{thm}
\label{41}
Let $f$ be a continuous function on a closed subset $\fF$ of $\C$. Then the functions
\bay
\label{do}
\d\mapsto\d^{-1}\O_f^\SA(\d)\quad\text{and}\quad\d\mapsto\d^{-1}\O_f^\CO(\d),\quad\d>0,
\ey
are nonincreasing. In particular,
$$
\O_f^\SA(\d_1\!+\!\d_2)\le\O_f^\SA(\d_1)+\O_f^\SA(\d_2)\!\quad
\text{and}\quad\!\O_f^\CO(\d_1\!+\!\d_2)\le\O_f^\CO(\d_1)+\O_f^\CO(\d_2),\!\quad\d_1,~\d_2>0.
$$
\end{thm}

%The proof is similar to the proof of Theorem 4.1 in \cite{AP4}

\Pf We consider first the case of the commutator modulus of continuity $\O_f^\CO$.
It suffices to verify that $\tau\O_{f}^{\CO}(\d/\tau)\le\O_{f}^{\CO}(\d)$
for $\d\in(0,\be)$ and $\tau\in(0,1)$.
It follows from Theorem \ref{com123}  that
\begin{align*}
\t\O_{f}^{\CO}(\d/\tau)&=
\t\sup\big\{\|f(N)R-Rf(N)\|:~\|R\|=1,~N\in\No(\fF),~\|NR-RN\|<\d/\tau\big\}\\[.2cm]
&=\sup\big\{\|f(N)R-Rf(N)\|:~\|R\|=\tau,~N\in\No(\fF),~\|NR-RN\|<\d\big\}\\[.2cm]
&\le\O_{f}^{\CO}(\d).
\end{align*}
The case of $\O_f^{\SA}$ can be treated in the same way if
we apply Theorem \ref{123} instead of Theorem \ref{com123}. $\bl$

\begin{cor}
The functions
$\O_f^\SA$ and $\O_f^\CO$ are continuous as functions from $(0,\be)$ to $[0,\be]$.
\end{cor}

\Pf It suffices to observe that if the function $h:(0,\be)\to[0,\be]$ is nondecreasing and the function
$\d\mapsto\d^{-1}h(\d)$ is nonincreasing, then $h$ is continuous. $\bl$

We can consider 3 more versions of commutator moduli of continuity. Let $f$ be a continuous function on a subset $\fF$ of $\C$. Put
\begin{align*}
&\O_f^{\Un}(\d)\df\sup\big\{\|f(N)U-Uf(N)\|:
~U\in\Un,~N\in\No(\fF),~\|NU-UN\|\le\d\big\};\\[.2cm]
&\O_f^{\Un\SA}(\d)\df\sup\big\{\|f(N)Q-Qf(N)\|:
~Q\in\USA,~N\in\No(\fF),~\|NQ-QN\|\le\d\big\};\\[.2cm]
&\O_f^{\Pro}(\d)\df\sup\big\{\|f(N)P-Pf(N)\|:
~P\in\Pro,~N\in\No(\fF),~\|NP-PN\|\le\d\big\}.
\end{align*}

\begin{thm}
\label{USA}
Let $f$ be a continuous function on a closed subset $\fF$ of $\C$. Then
$$
\O_f^{\Un}(\d)=\O_f^{\USA}(\d)=\O_f(\d)
\quad\mbox{and}\quad\O_f^{\Pro}(\d)=\frac12\O_f(2\d).
$$
\end{thm}

\Pf
Clearly, $\O_f^{\USA}\le\O_f^{\Un}\le\O_f$. The inequality $\O_f\le\O_f^{\USA}$ can be proved
in the same way as in the proof of the implication (iv)$\imp$(i) in Theorem \ref{sr}.
The equality $\O_f^{\Pro}(\d)=\frac12\O_f^{\USA}(2\d)$ follows from the fact that
$Q\in\USA$ if and only if $Q=2P-I$ for an orthogonal projection $P$.
$\bl$

Now we are going to show that the commutator moduli of continuity of $\O_{f}^{\SA}$
and $\O_{f}^{\CO}$ can be estimated in terms of operator Lipschitz norms
and commutator Lipschitz norms.

Recall that a subset $\L$ of $\C$ is called a {\it$\d$-net} for $\fF$ if
for every $z\in\fF$ there exists a $\l\in\L$ such that $|\l-z|\le\d$.

\begin{thm}
\label{kme+}
Let $f$ be a continuous function on a closed subset $\frak F$ of $\C$.
Suppose that $\fF_\d$ is a subset of $\fF$ that forms a $(\d/2)$-net of $\fF$. Then
$$
\O_{f}^{\SA}(\d)\le2\o_{f}(\d/2)+2\d\|f\|_{{\rm OL}(\frak F_\d)}
$$
and
$$
\O_{f}^{\CO}(\d)\le2\o_{f}(\d/2)+2\d\|f\|_{{\rm CL}(\frak F_\d)}.
$$
\end{thm}

The case $\fF\subset\R$ is Theorem 5.10 in \cite{AP6}. The general case can be proved
in the same way.

We need a lower estimate for the commutator modulus of continuity. The following theorem
can be considered as a version of Theorem 5.11 in \cite{AP2} for functions defined
on subsets of $\C$.

\begin{thm}
\label{kme-}
Let $f$ be a continuous function on a closed subset $\frak F$ of $\C$ and let $\d>0$.
Suppose that $\L$ and $\rm M$ are closed subsets of $\frak F$ such that
$(\L-{\rm M})\cap\bs{c}\d\,\dd\subset\{0\}$, where
$\bs{c}=\dfrac12\|\Psi\|_{\widehat L^1(\C)}$
and $\Psi$ is defined by {\em\rf{Phi}}.
Then
$$
\O_{f}^{\CO}(\d)\ge
\max\left\{\o_{f}(\d),\frac\d2\|\dg_0 f\|_{\fM_{\L,\rm M}}\right\}.
$$
\end{thm}

\Pf Clearly, $\o_{f}\le\O_{f}\le\O_{f}^{\CO}$. Note that
$$
\|\dg_0 f\|_{\fM_{\L,\rm M}}=
\sup_{a>0}\|\dg_0 f\|_{\fM_{\L\cap\clos(a\dd),{\rm M}\cap\clos(a\dd)}}.
$$
Thus it suffices to prove that
$$
\O_{f}^{\CO}(\d)\ge\frac\d2\|\dg_0 f\|_{\frak M_{\L,\rm M}}
$$
in the case when $\L$ and $\rm M$ are bounded.

Let $\e>0$. There exist Borel measures $\l$ on $\L$ and $\mu$
on $\rm M$, and a function $k$ in $L^2(\L\times{\rm M},\l\otimes\mu)$ such that
$$
\|k\|_{\mB_{\L,{\rm M}}^{\l,\mu}}=1,\quad
k\,\dg_0\in L^2(\L\times{\rm M},\l\otimes\mu),\quad\mbox{and}\quad
\|k\,\dg_0 f\|_{\mB_{\L,{\rm M}}^{\l,\mu}}\ge\|\dg_0 f\|_{\fM_{\L,{\rm M}}}-\e.
$$
We define the function $k_0$ in $L^2(\L\times{\rm M},\l\otimes\mu)$ by
$$
k_0(z,w)\df\left\{\begin{array}{ll}k(z,w),&\text {if}\,\,\,\,z\ne w,\\[.2cm]
0,&\text {if}\,\,\,\,z=w.
\end{array}\right.
$$
Then $k\dg_0 f=k_0\dg_0 f$ and by Corollary \ref{grunya},
$\|k_0\|_{\mB_{\L,{\rm M}}^{\l,\mu}}\le2$.
%Put $\Phi(x,y)\df f_\d(x-y)$ where $f_\d$ denotes the same as in Corollary \ref{Hxy}.
We define the normal operators $N_1:L^2(\L,\l)\to L^2(\L,\l)$ and
$N_2:L^2({\rm M},\mu)\to L^2({\rm M},\mu)$ by $(N_1f)(z)\df zf(z)$ and
$(N_2g)(w)\df wg(w)$.
Put
$$
h(z,w)\df\frac{1}{\bs{c}\d}\Psi\Big(\frac{z-w}{\bs{c}\d}\Big)k(z,w)=
\frac1{\bs{c}\d}\Psi\Big(\frac{z-w}{\bs{c}\d}\Big)k_0(z,w),
$$
where $\Psi$ is defined by \rf{Phi}.
Clearly,
$$
\|h\|_{\mB_{\L,{\rm M}}^{\l,\mu}}
\le\frac1{\bs{c}\d}
\left\|\Psi\Big(\frac{z-w}{\bs{c}\d}\Big)\right\|_{\fM_{\L,{\rm M}}^{\l,\mu}}
\|k\|_{\mB_{\L,{\rm M}}^{\l,\mu}}
\le\frac1{\bs{c}\d}\left\|\Psi\Big(\frac{z-w}{\bs{c}\d}\Big)\right\|_{\fM_{\C}}
=\frac{2}{\d}
$$
by Corollary \ref{corovz}.

Clearly, $N_1\mI_h-\mI_h N_2=\mI_{k_0}$ and
$f(N_1)\mI_h-\mI_h f(N_2)=\mI_{k_0\dg_0 f}$.
(Recall that $\mI_\f$ is the integral operator from $L^2({\rm M},\mu)$ into $L^2(\L,\l)$ with kernel
$\f\in L^2(\L\times{\rm M},\l\otimes\nu)$.)
Then
$$
\left\|\frac\d 2\mI_h\right\|=\frac\d 2\|h\|_{\mB_{\L,{\rm M}}^{\l,\mu}}\le1,
$$
$$
\left\|N_1\left(\frac\d 2\mI_h\right)-\left(\frac\d 2\mI_h\right) N_2\right\|=\frac\d 2\|k_0\|
_{\mB_{\L,{\rm M}}^{\l,\mu}}\le\d,
$$
%$$
%\|\mI_h\|=\|h\|_{\mB_{\L,{\rm M}}^{\l,\mu}}\le\frac C\d,
%$$
%$$
%\|\d \mI_h\|=\d \|h\|_{\mB_{\L,{\rm M}}^{\l,\mu}}\le1,
%$$
%$$
%\|N_1(\d \mI_h)-(\d \mI_h) N_2\|=\frac\d C\|k_0\|
%_{\mB_{\L,{\rm M}}^{\l,\mu}}\le\frac{2\d}C,
%$$
and
$$
\left\|f(N_1)\left(\frac\d 2\mI_h\right)-\left(\frac\d 2\mI_h\right) f(N_2)\right\|
=\frac\d 2\|k_0\dg_0 f\|_{\mB_{\L,{\rm M}}^{\l,\mu}}
\ge
\frac\d 2\big(\|\dg_0 f\|_{\fM_{\L,{\rm M}}}-\e\big).
$$
Hence, $\O_{f}^{\CO}(\d)\ge\frac\d 2\big(\|\dg_0 f\|_{\fM_{\L,{\rm M}}}-\e\big)$
for every $\e>0$.
$\bl$

Theorem \ref{kme-} allows us to obtain the following generalization of Theorem 4.17 in \cite{AP4}.

\begin{thm}
\label{JW}
Let $f$ be a continuous function on an unbounded closed subset $\fF$ of $\C$.
Suppose that $\O_{f}^\CO(\d)<\be$ for $\d>0$.
Then
the function $z\mapsto z^{-1}f(z)$ has finite limit as $|z|\to\be$, $z\in \fF$.
\end{thm}

\Pf We repeat the arguments of the proof of Theorem 5.12 in \cite{AP6}.
Assume the contrary. Then there exists a sequence $\{\l_n\}_{n=1}^\be$ in $\fF$ such that
$|\l_{n+1}|-|\l_n|>\bs{c}$ for all $n\ge1$, where $\bs{c}$ is the same as in
Theorem \ref{kme-}, $\lim_{n\to\be}|\l_n|=\be$, and the sequence $\{\l_n^{-1}f(\l_n)\}_{n=1}^\be$
has no finite limit. Put
$$
\L\df\{\l_n:~n\ge1\}.
$$
Then $\dg_0f\not\in\fM_{\L}$. This fact is contained implicitly in \cite{JW}.
Indeed, Theorem 4.1 in \cite{JW} implies that if $\dg_0f\in\fM_{\L}$, then
$f$ has complex derivative
at every nonisolated point of $\fF$. It can be shown that the same argument
gives us the differentiability
at $\be$ in the following sense: the function $z\mapsto z^{-1}f(z)$ has
finite limit as $|z|\to\be$,
provided the domain of $f$ is unbounded. Applying Theorem \ref{kme-}
for $M=\L$ and $\d=1$, we find that $\O_{f}^\CO(1)=\be$. $\bl$

\medskip

{\bf Remark.} Theorem \ref{JW} can also be proved in the same way as Theorem 4.17 in \cite{AP4}.

\medskip

Let $M$ and $N$ be (not necessarily bounded) normal operators
in a Hilbert space and let $R$ be a bounded operator on the same Hilbert space.
We say that the {\it operator $MR-RN$ is bounded} if
$R(\cd_N)\subset \cd_M$ and $\|MRu-RNu\|\le C\|u\|$ for every $u\in \cd_N$.
Then there exists a unique bounded operator $K$ such that
$Ku=MRu-RNu$ for all $u\in \cd_N$. In this case we write $K=MR-RN$.
Thus $MR-RN$ is bounded if and only if
\bay
\label{MN}
\big|(Ru,M^*v)-(Nu,R^*v)\big|\le C\|u\|\cdot\|v\|
\ey
for every $u\in \cd_N$ and $v\in \cd_{M^*}=\cd_{M}$. It is easy to see that
$MR-RN$ is bounded if and only if $N^*R^*-R^*M^*$ is bounded,
and $(MR-RN)^*=-(N^*R^*-R^*M^*)$.
In particular, we write $MR=RN$ if $R(\cd_N)\subset \cd_M$ and $MRu=RNu$ for every $u\in \cd_N$.
We say that $\|MR-RN\|=\be$ if $MR-RN$ is not a bounded operator.

Let us state now an analog of Lemma 4.4 in \cite{AP4} for normal operators. Recall that for a bounded operator $R$, the singular values $s_j(R)$ are defined by
$$
s_j(R)=\inf\big\{\|R-K\|:~\rank K\le j\big\},\quad j\ge0.
$$

\begin{lem}
\label{anbnrnorm}
Let $M$ and $N$ be (not necessarily bounded) normal operators and let
$R$ be an operator such that $\|R\|=1$. Then there exist
a sequence of operators $\{R_n\}_{n\ge1}$ and sequences of
bounded normal operators $\{M_n\}_{n\ge1}$ and $\{N_n\}_{n\ge1}$ such that

{\em(i)} $\s(M_n)\subset\s(M)\cup\{0\}$ and $\s(N_n)\subset\s(N)\cup\{0\}$;

{\em(ii)} the sequence $\{\|R_n\|\}_{n\ge1}$ is nondecreasing and
$$
\lim_{n\to\be}\|R_n\|=1;
$$

{\em(iii)}
$$
\lim_{n\to\be}R_n=R
$$
in the strong operator topology;

{\em(iv)} for every
continuous functions $f$ on $\C$, the sequence
$\big\{\big\|f(M_n)R_n-R_nf(N_n)\big\|\big\}_{n\ge1}$ is nondecreasing
and
$$
\lim_{n\to\be}\big\|f(M_n)R_n-R_nf(N_n)\big\|=\|f(M)R-Rf(N)\|;
$$

{\em(v)} if $f$ is a continuous function on $\C$ such that
$\|f(M)R-Rf(N)\|<\be$, then
$$
\lim_{n\to\be}f(M_n)R_n-R_nf(N_n)=f(M)R-Rf(N)
$$
in the strong operator topology;

{\em(vi)} if $f$ is a continuous function on $\C$ such that
$\|f(M)R-Rf(N)\|<\be$, then for every $j\ge0$, the sequence
$\big\{s_j\big(f(M_n)R_n-R_nf(N_n)\big)\big\}_{n\ge1}$ is nondecreasing and
$$
\lim_{n\to\be}s_j\big(f(M_n)R_n-R_nf(N_n)\big)=s_j\big(f(M)R-Rf(N)\big).
$$
\end{lem}

\Pf
Put $P_n\df E_M\big(n\dd)$, $Q_n\df E_N\big(n\dd)$  and $R_n\df P_nRQ_n$,
where $E_M$ and $E_N$ are the spectral measures of $M$ and $N$. Put $M_n\df P_nM=MP_n$ and
$N_n\df Q_nN=NQ_n$. Statements (i), (ii) end (iii) are evident.
Clearly,
$$
P_n\big(f(M_{n+1})R_{n+1}-R_{n+1}f(N_{n+1})\big)Q_n=f(M_n)P_nRQ_n-P_nRQ_nf(N_n),\quad n\ge1.
$$
Hence, the sequence $\big\{s_j\big(f(M_n)R_n-R_nf(N_n)\big)\big\}_{n\ge1}$ is nondecreasing
for every $j\ge0$. In particular, the sequence $\big\{\big\|f(M_n)R_n-R_nf(N_n)\big\|\big\}_{n\ge1}$ is nondecreasing.
Statement (v) follows from the identity
\bay
\label{pnqn}
P_n\big(f(M)R-Rf(N)\big)Q_n=f(M_n)P_nRQ_n-P_nRQ_nf(N_n),\quad n\ge1.
\ey
It is easy to see from \rf{pnqn} that (v)$\Rightarrow$(vi).

It remains to prove (iv). If $\lim_{n\to\be}\big\|f(M_n)R_n-R_nf(N_n)\big\|=\be$, then the result
follows from (v) with the help of an argument by contradiction.

Suppose that $\lim_{n\to\be}\big\|f(M_n)R_n-R_nf(N_n)\big\|<\be$.
Let $u\in\cd_{f(N)}$. Clearly, the sequence
$\{R_nf(N_n)u\}_{n=1}^\be=\{R_nf(N)u\}_{n=1}^\be$ converges.
Hence, the sequence \lb$\{f(M_n)R_nu\}_{n=1}^\be=\{f(M)R_nu\}_{n=1}^\be$
is bounded. Taking into account that $R_nu\to Ru$ as $n\to\be$, we find
that $Ru\in\cd_{f(M)}$ and $f(M_n)R_nu=f(M)R_nu\to f(M)Ru$ in the weak topology.
Thus, we have proved that $R(\cd_{f(N)})\subset \cd_{f(M)}$ and
$$
\|f(M)Ru-Rf(N)u\|\le\|u\|\lim_{n\to\be}\big\|f(M_n)R_n-R_nf(N_n)\big\|
$$
for all $u\in\cd_{f(N)}$. Hence,
$$
\|f(M)Ru-Rf(N)u\|\le\lim_{n\to\be}\big\|f(M_n)R_n-R_nf(N_n)\big\|.
$$
The opposite inequality follows from \rf{pnqn}. $\bl$

\begin{cor}
Theorem {\em\ref{sr}} remains valid for not necessarily bounded normal operators.
Moreover, each of the statements {\em (i)--(vii)} in Theorem {\em\ref{sr}}
is equivalent to the corresponding
statement for not necessarily bounded normal operators with spectrum in $\fF$.
\end{cor}

\Pf We can assume that $0\in\fF$. Note that statement (vii) for not necessarily bounded normal operators
with spectrum in $\fF$ implies the remaining statements. Thus it suffices to verify that statement (vii)
for bounded normal operators with spectrum in $\fF$ implies the same statement for arbitrary normal operators
with spectrum in $\fF$.
This immediately follows from Lemma \ref{anbnrnorm}. $\bl$

Lemma \ref{anbnrnorm} also implies that Theorems \ref{sr0}, \ref{123} and \ref{com123}
remain valid for not necessarily bounded normal operators..

Theorems \ref{123} and \ref{com123} for not necessarily bounded normal
operators
imply that we obtain the same commutator
moduli of continuity $\O_{f,\fF}^{\SA}$ and $\O_{f,\fF}^{\CO}$
if
we allow in the definitions of
$\O_{f,\fF}^{\SA}$
and $\O_{f,\fF}^{\CO}$ unbounded normal operators $N_1$ and $N_2$ .

We are not able to prove similar results for
$\O_{f,\fF}$, $\O_{f,\fF}^{\Un}$, $\O_{f,\fF}^{\Un\SA}$ and
$\O_{f,\fF}^{\Pro}$.

\begin{thm}
Let $\fF$ be an unbounded closed subset of $\C$ and let $f\in C(\fF)$.
Then
$$
\|f(N_1)-f(N_2)\|\le2\O_{f,\fF}(\|N_1-N_2\|)
$$
for arbitrary (not necessary bounded) normal
operators $N_1$ and $N_2$
with spectra in $\fF$.
\end{thm}
\Pf We have
$$
\|f(N_1)-f(N_2)\|\le\O_{f,\fF}^{\SA}(\|N_1-N_2\|)\le2\O_{f,\fF}(\|N_1-N_2\|).\quad\bl
$$

\medskip

{\bf Remark.} A similar result with the same constant 2 can be proved
for $\O_{f,\fF}^{\Un}$, $\O_{f,\fF}^{\Un\SA}$ and $\O_{f,\fF}^{\Pro}$.
Indeed, it
suffices to observe that Theorem \ref{USA} can be proved in the same way for
"unbounded" versions of
$\O_{f,\fF}$, $\O_{f,\fF}^{\Un}$, $\O_{f,\fF}^{\Un\SA}$ and
$\O_{f,\fF}^{\Pro}$.

\medskip

{\bf Definition 1.} For  a continuous function $f$ on a closed subset $\fF$ of $\C$, we consider the map
\bay
\label{fA}
N\mapsto f(N)
\ey
defined on the set of all (not necessarily bounded) normal operators with spectrum in $\fF$.
Let $N_0$ be a (not necessarily bounded) normal operator with spectrum in $\fF$.
We say that the mapping \rf{fA} is continuous at $N_0$ if
for arbitrary $\e>0$ there exists $\d>0$ such that
$\|f(N)-f(N_0)\|<\e$, whenever $N$ is a normal operator with spectrum in $\fF$
such that  $\|N-N_0\|<\d$.

We say that $f$ is {\it operator continuous}
if the map \rf{fA} is continuous at every (not necessarily bounded) normal operator $N$
with spectrum in $\fF$.

\medskip

It is easy to see that if $f$ is a continuous function on $\fF$,
then the map \rf{fA} is continuous at every bounded normal operator $N$ in $\No(\fF)$.
Indeed, this is obvious when $f$ is a polynomial of two real variables.
The result for arbitrary continuous functions follows from the
Stone--Weierstrass theorem applied to the polynomials on the closure of a bounded
neighborhood of $\s(N)$ in $\fF$.

\medskip

{\bf Definition 2.} Let $f$ be a continuous function on a closed subset $\fF$ of $\C$. It is called {\it uniformly operator continuous} if for every $\e>0$, there exists $\d>0$ such that
$\|f(N_1)-f(N_2)\|<\e$, whenever $N_1$ and $N_2$ are normal operators with spectra in $\fF$ such that $\|N_1-N_2\|<\d$.

\begin{thm}
\label{uoc}
Let $f$ be a bounded uniformly continuous function on $\C$. Then $f$ is uniformly operator continuous.
\end{thm}

\Pf Let $\o=\o_f$. Then $\o$ is a bounded modulus of continuity, and so $\o_*(\d)<\be$, $\d>0$.
The result follows now from Theorem 8.2 of \cite{APPS2}. $\bl$

\begin{cor}
Let $f$ be a bounded uniformly continuous function on a closed subset of $\C$. Then $f$ is uniformly operator continuous.
\end{cor}

\Pf It suffices to extend $f$
to a bounded and uniformly continuous function on $\C$.
To construct such an extension, one can use the first operator of continuation that was considered in \cite{S}, Ch. VI, \S 2.2. $\bl$

\begin{thm}
\label{omc0}
Let $f$ be an operator continuous function on a closed subset of $\C$. Then
$$
\lim_{\d\to0}\O_f(\d)=0,
$$
and so $f$ is uniformly operator continuous.
\end{thm}

The proof is similar to the proof of Theorem 8.2 in \cite{AP2}.

\

\section{\bf Estimates of commutator moduli of continuity}
\setcounter{equation}{0}
\label{ECMC}

\

In this section we obtain estimates of commutator moduli if continuity. In particular,
we show that for $\a\in(0,1)$, functions $f$ in the H\"older class $\L_\a(\R^2)$ must
be commutator H\"older of order $\a$, i.e.,
$\|f(N)R-Rf(N)\|\le\const\|Nr-RN\|^\a\|R\|^{1-\a}$ for every normal operator $N$ and
every bounded operator $R$.

%We start with several auxiliary results.

\begin{thm}
Let $f(z)=\ov z$ and $0<\d\le r$. Then
$$
C_1\d\log\frac{2r}\d\le\O_{f,\clos(r\dd)}^\CO(\d)\le C_2\d\log\frac{2r}\d,
$$
where positive numbers $C_1$ and $C_2$ are absolute constants.
\end{thm}
\Pf We first prove the upper estimate.
Note that the set $(\d/3)\mZ\cap\clos(r\dd)$ is a $(\d/2)$-net for $\clos(r\dd)$.
By Theorem \ref{kme+}  we have
$$
\O_{f,\clos(r\dd)}^{\CO}(\d)\le2\o_{f,\clos(r\dd)}(\d/2)+2\d\|f\|_{{\COL}(\frak F_\d)},
$$
where $\fF_\d=(\d/3)\mZ\cap\clos(r\dd)$. It remains to observe that
$2\o_{f,\clos(r\dd)}(\d/2)\le\d$
and to apply
Corollary \ref{log1}.
Now we prove the lower estimate. Put $\L\df{\rm M}\df\bs{c}\d\mZ\cap r\clos\dd$, where
$\bs{c}$ denotes the same as in Theorem \ref{kme-}.
By Theorems \ref{kme-} and \ref{niz},
we obtain
$$
\O_{f,\clos(r\dd)}^{\CO}(\d)\ge\frac\d2\|\dg_0 f\|_{\fM_{\L}}
\ge\const\d\log\frac{2{r}}{\bs{c}\d},
$$
provided ${r}>\bs{c}\d$. The case where $\d\le{r}<\bs{c}\d$ is evident. $\bl$

\begin{cor}
Let $N_1$ and $N_2$ be normal operators and let $R$ be a contraction.
Then
$$
\|RN_1^*-N_2^*R\|
\le\const\left(\log\frac{2(\|N_1\|+\|N_2\|)}{\|RN_1-N_2R\|}\right)\|RN_1-N_2R\|.
$$
\end{cor}

\Pf Put $r\df\|N_1\|+\|N_2\|$. It remains to apply Theorem \ref{com123}. $\bl$

We use the following notation:
$$
\D_X\df\{(x,x):\,\,x\in X\}
$$
for a set $X$.

\begin{lem}
\label{part}
For every $a>0$, there exists a Borel partition $\big\{G_a^{(j)}\big\}_{j=0}^9$ of $\C^2$
such that
$$
\Big\{(z,w)\in\C^2:\dist((z,w),\D_\C)<a\Big\}\subset
\bigcup_{j=1}^9G_{a}^{(j)}\subset\Big\{(z,w)\in\C^2:\dist((z,w),\D_\C)<3a\Big\},
$$
$\big\|\chi_{G_a^{(j)}}\big\|_{\fM_{\C,\C}}=1$ for $j=1,2,\dots,9$,
and
$\big\|\chi_{G_a^{(0)}}\big\|_{\fM(\C\times\C)}\le10$.
\end{lem}

\Pf We consider first the standard partition of $\C$ in
the squares $Q_\fn$ with side of length $1$:
$$
\C=\bigcup_{\fn\in\mZ}Q_\fn,
$$
where
$$
Q_0\df\{z\in\C:0\le\re z,\,\im z<1\}\quad\text{and}\quad Q_\fn\df\fn+Q_0.
$$
Denote by $\mathscr Q$  the set of all squares
$Q_\fn$ with $\fn\in\mZ$.
Given $\fn\in\mZ$, we define the set $\X_\fn$ by $\X_\fn\df\bigcup\limits_{Q\in\mathscr Q}\big(Q\times(\fn+Q)\big)$.
Clearly, the family $\{\X_\fn\}_{\fn\in\mZ}$ forms a partition of $\C^2$,
and $\big\|\chi_{_{\!\X_\fn}}\big\|_{\fM_\C}=1$ for every $\fn\in\mZ$ by \rf{1mult}.
We enumerate the sets $\mX_\fn$ with $|\fn|\le\sqrt2$ by a sequence
$\big\{G^{(j)}\big\}_{j=1}^9$ and put $G^{(0)}\df\C^2\setminus\bigcup_{j=1}^9 G^{(j)}$.
Then $\big\|\chi_{G_a^{(j)}}\big\|_{\fM_\C}=1$ for $j=1,2,\dots,9$ and
$$
\big\|\chi_{G_a^{(0)}}\big\|_{\fM_\C}\le1+\sum_{j=1}^9\big\|\chi_{G_a^{(j)}}\big\|_{\fM_\C}=10.
$$
It is easy to see that
$$
\big\{(z,w)\in\C^2:|z-w|\le1\big\}\subset\bigcup_{j=1}^9G_{a}^{(j)}\subset\big\{(z,w)\in\C^2:|z-w|<2\sqrt2\big\}
$$
and $\sqrt2\,\,\dist\big((z,w),\D_\C\big)=|z-w|$. Putting $G_a^{(j)}\df G^{(j)}$
for $j=0,1,\dots9$,
we obtain the desired result for $a=\frac1{\sqrt2}$.
The general case can be reduced to this special case with the help of dilations. $\bl$

\begin{lem}
\label{59}
Suppose that $\s$ and $\d$ are positive numbers and let $f$ be a bounded continuous function on $\C$ such that
$$
\supp\F f\subset\{\z\in\C:~|\z|\le\s\}.
$$
Then
$$
\|f\|_{\COL(\d\mZ)}\le\|\dg_0 f\|_{\fM_{\,\d\mZ}}\le\const\s\log\left(1+\frac1{\s \d}\right)\|f\|_{L^\be}.
$$
\end{lem}

\Pf The case $\s \d\ge1$ reduces to Corollary \ref{lat}. Suppose that $\s \d\le1$.
Clearly, it suffices to consider the case when $\d=1$.
It follows from Theorem 5.1 in \cite{APPS2} that there exist functions
$g_1,\,g_2\in\fM_\C$ such that
$f(z)-f(w)=g_1(z,w)(z-w)+g_2(z,w)(\ov z-\ov w)$ and $\|g_1\|_{\fM_\C}
+\|g_2\|_{\fM_\C}\le\const\s\|f\|_{L^\be}$. We do not need the continuity of $g_1$ and $g_2$.
Thus we may assume that
$g_1(z,z)=g_2(z,z)=0$ for every $z\in\C$. Clearly,
$$
(\dg_0 f)(z,w)=g_1(z,w)+(\ov z-\ov w)\l(z-w)g_2(z,w).
$$
%Thus it suffices to estimate $\big\|(\ov z-\ov w)\l(z-w)g_2(z,w)\big\|_{\fM_\mZ}$.

Let $G_a^{(j)}$ denote the same as in Lemma \ref{part} with $a=\s^{-1}$. Then
\begin{align*}
(\dg_0 f)(z,w)=&\chi_{G_a^{(0)}}(z,w)(\dg_0 f)(z,w)+\sum_{j=1}^9\chi_{G_a^{(j)}}(z,w)g_1(z,w)\\[.2cm]
&+\sum_{j=1}^9(\ov z-\ov w)\l(z-w)\chi_{G_a^{(j)}}(z,w)g_2(z,w).
\end{align*}
We have
$$
\chi_{G_a^{(0)}}(z,w)(\dg_0 f)(z,w)=f(z)\chi_{G_a^{(0)}}(z,w)\l(z-w)-f(w)\chi_{G_a^{(0)}}(z,w)\l(z-w).
$$
Hence,
\begin{align*}
\big\|\chi_{G_a^{(0)}}(z,w)(\dg_0 f)(z,w)\big\|_{\fM_\mZ}&\le2\|f\|_{L^\be}
\big\|\chi_{G_a^{(0)}}(z,w)\l(z-w)\big\|_{\fM_\mZ}\\[.2cm]
&=\frac{\sqrt2}a\|f\|_{L^\be}\left\|\chi_{G_a^{(0)}}(z,w)\Phi\left(\frac{z-w}{\sqrt2 a}\right)\right\|_{\fM_\mZ}
\le\const\s\|f\|_{L^\be}
\end{align*}
%Observe that
%$$
%\|\chi_{G_a^{(0)}}(z,w)l(z-w)\|_{\fM(\mZ\times\mZ)}\le Ca^{-1}=C\s
%$$
by Corollary \ref{corovz}. Let $j\ge1$. Applying Lemma \ref{part}, we obtain
$$
\big\|\chi_{G_a^{(j)}}(z,w)g_1(z,w)\big\|_{\fM_\mZ}\le\|g_1(z,w)\|_{\fM_\C}
\le\const\s\|f\|_{L^\be}.
$$
It remains to estimate $\big\|(\ov z-\ov w)\l(z-w)\chi_{G_a^{(j)}}(z,w)g_2(z,w)\big\|_{\fM_\mZ}$
for $j\ge1$.
Lemma \ref{log0} implies the following inequality
$$
\big\|(\ov z-\ov w)\l(z-w)\chi_{G_a^{(j)}}(z,w)\big\|_{\fM_\mZ}\le\const\log(1+\s^{-1}).
$$
Hence,
\begin{align*}
\big\|(\ov z-\ov w)\l(z-w)\chi_{G_a^{(j)}}(z,w)g_2(z,w)\big\|_{\fM_\mZ}
&\le\const\log(1+\s^{-1})\|g_2(z,w)\|_{\fM_\C}\\[.2cm]
&\le\const\s\log(1+\s^{-1})\|f\|_{L^\be}.\quad\bl
\end{align*}

We need the following well-known result in approximation theory
(see e.g., Theorem 2.1 in \cite{APPS2}):

\begin{thm}
\label{prib}
For every uniformly continuous
function $f$ on $\C$ and every positive $d$ there exists a uniformly continuous
function $f_d$ such that
$$
\supp\F f_d\subset\{\z\in\C:~|\z|\le d^{-1}\},\quad\o_{f_d}\le C\o_f,\quad
\mbox{and}\quad\|f-f_d\|_{L^\be}\le C\o_f(d),
$$
where $C$ is an absolute constant.
\end{thm}

Recall that for a modulus of continuity $\o$,
$$
\o_{*}(\d)\df\d\int_\d^\be\frac{\o(t)}{t^2}\,dt= \int_1^\be\frac{\o(\d t)}{t^2}\,dt
$$
and
$$
\o_{**}(\d)\df(\o_{*})_*(\d)=\delta\int_\delta^\be\frac{\o(t)\log(t/\d)}{t^2}\,dt=
\int_1^\be\frac{\o(\d t)\log t}{t^2}\,dt.
$$

Clearly, for nondecreasing $\o$, we have $\o\le\o_*\le\o_{**}$.

\begin{thm}
\label{Lip}
Let $f$ be a uniformly continuous function on $\C$ such that
\lb$(\o_f)_{**}(\d)<\be$ for $\d>0$. Then $\dg_0 f\in\fM_{\,\d\mZ}$ and
$$
\|\dg_0 f\|_{\fM_{\,\d\mZ}}\le \frac{C(\o_f)_{**}(\d)}\d
$$
for every $\d>0$, where $C$ is an absolute constant.
\end{thm}

\Pf By Theorem \ref{prib}, there exists a sequence of uniformly continuous functions $\{f_n\}_{n=0}^\be$
such that
$$
\supp\F f_n\subset\big\{\z\in\C:~|\z|\le 2^{-n}\d^{-1}\big\}, \quad
\o_{f_n}\le C\o_f,\quad\mbox{and}\quad
\|f-f_n\|_{L^\be}\le C\o_f(2^n\d),
$$
where $C$ is an absolute constant.

We may assume that $f_n(0)=f(0)$.
Then we can find a subsequence $\{f_{n_j}\}_{j=0}^\be$ that
converges everywhere on $\C$. Note that if $g(z)=\lim\limits_{j\to\be}f_{n_j}(z)$ for all
$z\in\C$, then $\o_{g}\le C\o_f$ and $\supp\F g=\{0\}$. Hence, $g$ is a constant
function because $(\o_f)_{**}(\d)<\be$.
This implies that $\lim_{n\to\be}\dg_0 f_n(z,w)=0$ for all $(z,w)\in\C^2$.
Consequently,
$$
\|\dg_0 f\|_{\fM_{\,\d\mZ}}\le\|\dg_0(f-f_0)\|_{\fM_{\,\d\mZ}}
+\sum_{n=0}^\be\|\dg_0(f_n-f_{n+1})\|_{\fM_{\,\d\mZ}}.
$$
By Corollary \ref{lat},
$$
\|\dg_0(f-f_0)\|_{\fM_{\,\d\mZ}}\le\const\frac{\o_f(\d)}\d.
$$
By Lemma \ref{59},
$$
\|\dg_0(f_n-f_{n+1})\|_{\fM_{\,\d\mZ}}\le\const2^{-n}\d^{-1}(n+1)\o_f(2^n\d).
$$
It remains to observe that
\begin{align*}
\o_f(\d)+\sum_{n=0}^\be 2^{-n}(n+1)\o_f(2^n \d)
&\le\const\int_1^\be\frac{\o_f(t\d)\log(t+1)}{t^2}dt\\[.2cm]
&\le\const\o_f(\d)+\const\int_2^\be\frac{\o_f(t\d)\log t}{t^2}dt\\[.2cm]
&\le\const\int_1^\be\frac{\o_f(t\d)\log t}{t^2}dt=\const(\o_f)_{**}(\d).\quad \bl
\end{align*}
\begin{cor}
\label{lip1}
Let $\d>0$ and let $f$ be a uniformly continuous function on $\C$ such that
$(\o_f)_{**}(\d)<\be$, $\d>0$. Suppose that $N_1$ and $N_2$ are
normal operators  with spectra
in $\d\mZ$ and $R$ is a contraction such that the operator $RN_1-N_2R$ is bounded. Then
$$
\|Rf(N_1)-f(N_2)R\|\le\const \d^{-1}(\o_f)_{**}(\d)\|RN_1-N_2R\|.
$$
\end{cor}

\Pf The result follows from Theorem \ref{Lip}. $\bl$
%and from the first inequality in \rf{nva}. $\bl$

\begin{thm}
\label{513}
Let $f$ be a uniformly continuous function on $\C$.
Then $\O_f^{\CO}\le\const{\o_f}_{**}$.
\end{thm}

\Pf It suffices to apply Theorem \ref{Lip} and Theorem \ref{kme+} for $\fF=\C$
and $\fF_\d=(2\d/3)\mZ$. $\bl$

\begin{thm}
\label{o**}
Let $f$ be a uniformly continuous function on $\C$ such that
$(\o_f)_{**}(\d)<\be$, $\d>0$. If $N_1$ and $N_2$ are normal operators and $R$
is a contraction such that the operator $RN_1-N_2R$ is bounded, then
$$
\|Rf(N_1)-f(N_2)R\|\le\const\,(\o_f)_{**}\big(\|RN_1-N_2R\|\big).
$$
\end{thm}

\Pf The result follows from Theorem \ref{513} and Theorem \ref{com123}.  $\bl$

\begin{cor}
Let $\o$ be a modulus of continuity such that $\o_*\le C\o$ for a positive number $C$
and let $f\in\L_\o(\R^2)$. If $N_1$ and $N_2$ are normal operators and $R$
is a contraction such that the operator $RN_1-N_2R$ is bounded, then
$$
\|Rf(N_1)-f(N_2)R\|\le\const\o(\|RN_1-N_2R\|).
$$
\end{cor}

\Pf It suffices to observe that $\o_{**}\le C^2\o$ and  apply Theorem \ref{o**}. $\bl$

\begin{thm}
Let $0<\a<1$ and $f\in\L_\a(\R^2)$. Then
$$
\|Rf(N_1)-f(N_2)R\|\le\const\frac1{(1-\a)^2}\|RN_1-N_2R\|^\a\|R\|^{1-\a},
$$
whenever $N_1$ and $N_2$ are normal operators  and $R$ is bounded operator such that
the operator $RN_1-N_2R$ us bounded.
\end{thm}

\Pf The result can be easily deduced from Theorem \ref{o**}. $\bl$

\begin{cor}
Let $f\in{\rm Lip}(\C)\cap L^\be(\C)$. Then
$$
\|Rf(N_1)-f(N_2)R\|
\le\const\|f\|_{\rm Lip}\|RN_1-N_2R\|\cdot
\log^2\left(2+\dfrac{\|f\|_{L^\be}}{\|f\|_{\rm Lip}\|RN_1-N_2R\|}\right),
$$
whenever $N_1$ and $N_2$ are normal operators  and $R$ is bounded operator such that
the operator $RN_1-N_2R$ us bounded.
\end{cor}

\Pf Again, the result immediately follows from Theorem \ref{o**}. $\bl$

\medskip

{\bf Remark.} It is interesting to compare above
results with the results of \cite{APPS2}
quoted in the introduction that estimate the quasicommutator norms $\|Rf(N_1)-f(N_2)R\|$
in terms of $\max\big\{\|RN_1-N_2R\|,\|RN^*_1-N_2^*R\|\big\}$.

\medskip

Theorem \ref{513} and all results that follows from this theorem
can be generalized, in the spirit of Theorem \ref{W}, to the  case of functions defined on
a closed subset of $\C$.
We state only a version of Theorem \ref{513}.
\begin{thm}
Let $f\in\L_\o(\fF)$, where $\fF$ is a closed subset of $\C$
and $\o$ is a modulus of continuity. Then
$$
\O_{f,\fF}^{\CO}\le C\|f\|_{\L_\o(\fF)}\o_{**},
$$
where $C$ is an absolute constant.
\end{thm}

The proof repeats the proof of Theorem \ref{W}.

\

\section{\bf Constants in operator H\"older inequalities}
\setcounter{equation}{0}
\label{alpha}

\

As we have mentioned in the the introduction, it was proved in \cite{APPS2} that the
H\"older class $\L_\a(\R^2)$, $0<\a<1$, coincides with the class of operator H\"older
functions of order $\a$. Moreover, the best constant in inequality \rf{a^-1} can be
estimated from above in terms of $\const(1-\a)^{-1}$.
In this section we obtain a lower estimate for the constant
in inequality \rf{a^-1}.

Consider the operator H\"older semi-norm on $\L_\a(\R^2)$:
$$
\|f\|_{\rm O\L_\a}\df\sup\left\{\frac{\|f(N_1)-f(N_2)\|}{\|N_1-N_2\|^\a}\right\},
\quad0<\a<1,
$$
where the supremum is taken over all bounded normal operators $N_1$ and $N_2$
such that $N_1\ne N_2$. Denote by ${\rm O}\L_\a(\R^2)$ the space $\L_\a(\R^2)$ equipped with the semi-norm
$\|\cdot\|_{{\rm O}\L_\a}$.

Let ${\frak h}_\a$ be the norm of the identity operator from
$\L_\a(\R^2)$ to ${\rm O}\L_\a(\R^2)$, i.e.,
$$
{\frak h}_\a\df\sup\{\|f\|_{\rm O\L_\a}:~f\in\L_\a(\R),~\|f\|_{\L_\a}\le1\}.
$$
Recall that it was proved in \cite{APPS2} that ${\frak h}_\a\le\const(1-\a)^{-1}$.

\begin{thm}
There exists a positive constant $C$ such that
${\frak h}_\a\ge C(1-\a)^{-1/2}$ for all $\a\in(0,1)$.
\end{thm}

\Pf We prove the desired low estimate even in the case of
functions on $\R$ and self-adjoint operators.

By Theorem 9.9 in \cite{AP6}, there exists a function $f\in C^\be(\R)$ such that
$\|f\|_{L^\be}\le1$, $\|f^\prime\|_{L^\be}\le1$ and $\O_f(\d)\ge c\,\d\sqrt{\log\frac2\d}$
for every $\d\in(0,1)$, where $c$ is a positive constant. Clearly,
$$
\frac{|f(x)-f(y)|}{|x-y|^\a}=\frac{|f(x)-f(y)|^\a}{|x-y|^\a}|f(x)-f(y)|^{1-\a}\le2^{1-\a},
\quad x\ne y,
$$
and so
$\|f\|_{\L_\a}\le2^{1-\a}$ for $\a\in(0,1)$. Hence,
$$
{\frak h}_\a\ge\frac c{2^{1-\a}}\d^{1-\a}\sqrt{\log\frac2\d}
$$
for arbitrary $\a$ and $\d$ in $(0,1)$. Substituting $\d=2\exp(-\frac1{1-\a})$, we obtain
$$
{\frak h}_\a\ge\frac{c}e(1-\a)^{-1/2}
$$
for $\a\in(0,1)$. $\bl$

\

\

\footnotesize
\noindent
\begin{tabular}{p{8.7cm}p{5cm}}
A.B. Aleksandrov & V.V. Peller \\
St.Petersburg Branch & Department of Mathematics  \\
Steklov Institute of Mathematics  & Michigan State University\\
Fontanka 27  & East Lansing \\
191023 St-Petersburg &Michigan 48824\\
Russia & USA
\end{tabular}

\end{document}